\documentclass[a4paper,10pt]{amsart}
\usepackage{egen, egenengelsk}
\usepackage{hyperref}

\usepackage{tikz-cd}

\newcommand{\fa}{\mathfrak{a}}
\renewcommand{\wt}{\widetilde}
\newcommand{\std}{\mathrm{std}}
\newcommand{\Pf}{\mathrm{Pf}}
\newcommand{\op}{\mathrm{op}}
\renewcommand{\ss}{\mathrm{ss}}
\author{Jørgen Vold Rennemo}
\address{All Souls College, Oxford, OX1 4AL, UK}
\email{jvrennemo@gmail.com}

\title[The fundamental theorem of HPD via VGIT]{The fundamental theorem of homological projective duality via variation of GIT stability}
\begin{document}
\begin{abstract}
We reprove Kuznetsov's ``fundamental theorem of homological projective duality'' using LG models and variation of GIT stability.
This extends the validity of the theorem from smooth varieties to nice subcategories of smooth quotient stacks, and moreover shows that the line bundles polarising the varieties in HP duality do not have to be globally generated.

The proof combines ideas from Kuznetsov's original proof with ideas from work of Ballard--Deliu--Favero--Isik--Katzarkov.
\end{abstract}

\maketitle
\section{Introduction}
In 2005, Kuznetsov formulated the concept of ``Homological projective duality'' \cite{kuznetsov_homological_2007}.
The starting data for the theory is a pair of varieties $X$ and $Y$ which map to dual projective spaces $\PP V$ and $\PP V^{\vee}$, and which are equipped with special semi-orthogonal decompositions of their derived categories called ``Lefschetz decompositions''.
With this data one can define the notion of $X$ and $Y$ being HP dual, which is a relation between their derived categories.

Choosing a linear subspace $L \subseteq V$ gives the base changed varieties $Y_{L} = Y|_{\PP L}$ and $X_{L^{\perp}} = X|_{\PP L^{\perp}}$.
The ``fundamental theorem of HP duality'' states that if $X$ and $Y$ are HP dual and the base changed varieties have the expected dimension, then there exist explicit semiorthogonal decompositions of $D^{b}(Y_{L})$ and $D^{b}(X_{L^{\perp}})$, such that the two decompositions have a piece $\cC_{L}$ in common.

This result underlies many examples of interesting semi-orthogonal decompositions in geometry -- we refer to \cite{kuznetsov_homological_2007, kuznetsov_semiorthogonal_2014, rennemo_homological_2015, thomas_notes_2015} for further explanations of and examples in the theory of HP duality.

In this paper, we give a new proof of the fundamental theorem of HP duality, applying the technology of variation of GIT stability for LG models.
The proof is inspired by the work of Ballard--Deliu--Favero--Isik--Katzarkov \cite{ballard_homological_2013}.
In Section \ref{sec:relation} we explain the relation between that paper and this one further.

\subsection{Results}
We'll prove the theorems of HP duality for an admissible subcategory $\cD \subseteq D^{b}(X)$, where $X$ is a smooth quotient stack.
To be precise, we assume $X = [A/G]$, where $A$ is a smooth, quasi-projective variety and $G$ is a reductive group.
The case treated in \cite{kuznetsov_homological_2007} is the one where $\cD = D^{b}(X)$ and $X$ is a smooth, projective variety.
The extra generality is needed in applications to the HP duality results of \cite{rennemo_homological_2015, rennemo_hori-mological_2016}, where we work with non-commutative resolutions of singular varieties.
Apart from this generalisation and some weakening of hypotheses, discussed in Section \ref{sec:relation}, our result is \cite[Thm.~1.1]{kuznetsov_homological_2007}.

We assume that $X$ is ext finite, i.e.~that for coherent sheaves $\cE,\cF$ on $X$, the space $\oplus_{i\in\ZZ}\Ext^{i}(\cE,\cF)$ is finite-dimensional.
Choose a line bundle $\cL$ on $X$.
Let $\cD \subseteq D^{b}(\cX)$ be a subcategory equipped with a Lefschetz decomposition, i.e.~a sequence of admissible subcategories $\cA_{0} \supseteq \cA_{1} \supseteq \cdots \supseteq \cA_{n}$ giving a semiorthogonal decomposition
\[
\cD = \langle \cA_{0}, \cA_{1}(\cL), \ldots, \cA_{n}(n\cL) \rangle.
\]
We assume that $\cD$ is closed under tensoring with $\cL^{\pm 1}$, and that the $\cA_{i}$ are saturated -- these conditions hold in particular when $\cD = D^{b}(X)$ with $X$ a smooth, projective variety.

Let $V = H^{0}(X,\cL)^{\vee}$, choose a linear subspace $L \subset V^{\vee}$, and let $X_{L^{\perp}}$ be the base locus of the linear system $\PP L \subset |\cL|$.
We'll define a category $\cD_{L^{\perp}}$ which functions as a ``categorical base change'' of $\cD$ to $X_{L^{\perp}}$.
The actual definition is phrased in terms of matrix factorisation categories, but for an informal picture of $\cD_{L^{\perp}}$, it should be thought of as a subcategory of $D^{b}(X_{L^{\perp}})$, where $X_{L^{\perp}}$ is considered as a derived stack.

We'll also define a category $\cD^{\vee}$, which we call the HP dual category of $\cD$.
The category $\cD^{\vee}$ naturally lives over $\PP V^{\vee}$, and we write $\cD^{\vee}_{L}$ for its base change to $\PP L$.
These categories $\cD_{L^{\perp}}, \cD^{\vee}$ and $\cD^{\vee}_{L}$ are all defined as subcategories of certain matrix factorisation categories; we postpone their precise definition to Section \ref{sec:HPDSection}.

Let $\cM = \cO_{\PP V^{\vee}}(1)$, let $l = \dim L$, and let $m = \min\{i \ge 0 \mid \cA_{i} \not= \cA_{0} \}$.
The next result is the fundamental theorem of HP duality.
Let $l = \dim L$.
\begin{nthm}
\label{thm:fundamentalTheoremOfHPD}
There exists admissible subcategories $\cB_{m} \subseteq \cdots \subseteq \cB_{l-1} \subseteq \cD_{L}^{\vee}$ and a triangulated category $\cC_{L}$ with fully faithful inclusions $\cC_{L} \to \cD_{L^{\perp}}$ and $\cC_{L} \to \cD^{\vee}_{L}$, such that there are semiorthogonal decompositions
\[
\cD_{L^{\perp}} = \langle \cC_{L}, \cA_{l}(l\cL), \ldots, \cA_{n}(n\cL) \rangle
\]
and
\[
\cD^{\vee}_{L} = \langle \cB_{l-1}(-(l-1)\cM), \ldots, \cB_{m}(-m\cM), \cC_{L} \rangle.
\]
\end{nthm}
The $\cB_{i}$ are determined by the $\cA_{i}$ via an equivalence $\cB_{i} \cong \cA_{0} \cap \cA_{i}^{\perp}$.
Note that in contrast to \cite[Thm.~1.1]{kuznetsov_homological_2007} we don't require $X_{L^{\perp}}$ to have the expected dimension.
This is essentially because we've defined our categories $\cD_{L^{\perp}}$ and $\cD^{\vee}_{L}$ in such a way that base change is always well behaved.

Let $v = \dim V$.
Note that if $\cC_{V^{\vee}} = 0$, which holds in particular if $\cL$ is globally generated, the $\cB_{i}$ give a Lefschetz decomposition of the dual
\[
\cD^{\vee} = \langle \cB_{v-1}(-(v-1) \cM), \ldots, \cB_{m}(-m\cM) \rangle.
\]

In the original formulation of the fundamental theorem \cite[Thm.~1.1]{kuznetsov_homological_2007}, one takes $\cD = D^{b}(X)$, and assumes that there exists a variety $Y$ over $\PP V^{\vee}$ such that there is a $\PP V^{\vee}$-linear equivalence $D^{b}(Y) \cong \cD^{\vee}$.
The theorem then states the relation of Thm.~\ref{thm:fundamentalTheoremOfHPD}, with $\cD_{L^{\perp}}$ and $\cD_{L}^{\vee}$ replaced by $D^{b}(X_{L^{\perp}})$ and $D^{b}(Y_{L})$, under the additional assumption that $X_{L^{\perp}}$ and $Y_{L}$ have the expected dimensions.
It follows from Prop.~\ref{thm:KnorrerPeriodicity} that $D^{b}(X_{L^{\perp}}) \cong \cD_{L^{\perp}}$ if $X_{L}$ has the expected dimension.
Hence in order to recover the geometric formulation of the fundamental theorem, we need to see that $\cD^{\vee}_{L} \cong D^{b}(Y_{L})$ when $Y_{L}$ has expected dimension.
We give an argument for this based on results from \cite{kuznetsov_base_2011} in Section \ref{sec:baseChange}.

\subsection{Outline of the argument}
Let us in this outline restrict to the special case where $\cD = D^{b}(X)$, since the general argument differs only by taking appropriate subcategories of all categories involved.

Let $\cZ_{+}$ be the total space of the bundle $L \otimes \cL^{\vee}$ over $X$.
There is a natural potential $W \colon \cZ_{+} \to \CC$, and by Kn\"orrer periodicity (Prop.~\ref{thm:KnorrerPeriodicity}), the matrix factorisation category $D^{b}(\cZ_{+},W) \cong D^{b}(X_{L^{\perp}})$ if $X_{L^{\perp}}$ has the expected dimension.
Motivated by this fact, for general $L$ we \emph{define} $\cD_{L^{\perp}}$ to be $D^{b}(\cZ_{+},W)$

Next we note that $\cZ_{+}$ is a GIT quotient for a $\CC^{*}$-quotient problem.
Namely, let $\wt X$ be the total space of $\cL^{\vee}$, and let $Z = \wt X \times L$.
Let $\CC^{*}$ act on $\wt X \times L$ by scaling the fibres of $\wt X \to X$ with weight 1, and by scaling $L$ with weight $-1$.
We may then consider the quotient stack $\cZ = [Z/\CC^{*}]$, and can view $\cZ_{+} \subset \cZ$ as the open substack of semistable points for one choice of GIT stability.
Taking another choice of stability gives another open substack $\cZ_{-} = [\wt X \times (L \setminus 0)/\CC^{*}] \subseteq \cZ$.

The stack $\cZ_{-}$ maps to $\PP L$, and the fibre at $[f] \in \PP L$ is isomorphic to the LG model $(\wt X, f)$, where $f \in H^{0}(X,\cL)$ naturally defines a function on $\wt X$.
Applying Kn\"orrer periodicity, this shows that $D^{b}(\cZ_{-},W)$ is equivalent to $D^{b}(\cH)$, where $\cH \stackrel{i}{\into} X \times \PP L$ is the incidence stack of pairs $(x,[f])$ with $f(x) = 0$.
Kuznetsov's definition of the HP dual identifies the base change of the HP dual to $\PP L$ with a subcategory of $\cH$.
It can be described as the subcategory of those objects $\cE$ such that $(i_{*}\cE)|_{X \times [f]} \in \cA_{0}$ for all $[f] \in \PP L$.
Motivated by this, we define $\cD_{L}^{\vee}$ as the corresponding subcategory of $D^{b}(\cZ_{-},W)$.

The goal is now to use the GIT presentation to relate $\cD_{L^{\perp}}$ and $\cD_{L}^{\vee}$.
Applying the ideas of \cite{ballard_variation_2012, halpern-leistner_derived_2015, segal_equivalence_2011}, we find fully faithful inclusions $D^{b}(\cZ_{\pm},W) \into D^{b}(\cZ,W)$.
We then let $\cC_{L}$ be the intersection of the images of $\cD_{L^{\perp}}$ and $\cD_{L}^{\vee}$ inside $D^{b}(\cZ,W)$, which gives the inclusions of $\cC_{L}$ from Theorem \ref{thm:fundamentalTheoremOfHPD}.
The main technical difficulty in this argument is that the inclusion $D^{b}(\cZ_{+},W) \into D^{b}(\cZ,W)$ coming from the general results of \cite{ballard_variation_2012, halpern-leistner_derived_2015, segal_equivalence_2011} doesn't give the right intersection of the image categories, and so must be modified.
The correct inclusion is obtained in Section \ref{sec:windows}.

Obtaining the rest of Theorem \ref{thm:fundamentalTheoremOfHPD}, i.e.~describing the semi-orthogonal complements to $\cC_{L}$, is a matter of analysing the subcategories of $D^{b}(\cZ,W)$ in more detail.

\subsection{Relation to existing work}
\label{sec:relation}
This approach to HP duality is essentially a combination of techniques from Kuznetsov's original paper \cite{kuznetsov_homological_2007} and Ballard--Deliu--Favero--Isik--Katzarkov's paper \cite{ballard_homological_2013}.
Ballard et al.~show, using VGIT for LG models, that given a nice GIT quotient $X = [A^{\ss}/G]$, one can construct a Lefschetz decomposition of $D^{b}(X)$ and describe its HP dual category as a certain matrix factorisation category.
Their technique gives an independent proof of the fundamental theorem of HP duality for such $X$ with the constructed Lefschetz decomposition, i.e.~they don't rely on the results of \cite{kuznetsov_homological_2007}.
In contrast, we say nothing about constructing Lefschetz decompositions, but instead consider an arbitrary Lefschetz decomposition to be part of the starting data.

In \cite{kuznetsov_homological_2007}, Kuznetsov presents the definitions and framework of HP duality and proves our Thm.~\ref{thm:fundamentalTheoremOfHPD} in the case where $\cD = D^{b}(X)$ with $X$ a smooth, projective variety and assuming there exists a variety $Y$ over $\PP V^{\vee}$ such that $D^{b}(Y) \cong \cD^{\vee}$.
A more general ``categorical'' version of HP duality was worked out in Kuznetsov's habilitation thesis, where the assumption that such a $Y$ exists was dropped.

We follow \cite{kuznetsov_homological_2007} in definitions and parts of the proof, and several constructions are translated from that paper into the language of factorisation categories.
However, our proof does make it convenient to diverge in two points of notation, so let's list these.
Firstly, our categories $\cB_{i}$ are not the same as the ones used by Kuznetsov, because in the terminology of Section \ref{sec:dualLefschetz}, our $\cB_{i}$ form a right Lefschetz decomposition of $\cD^{\vee}$, whereas the $\cB_{i}$ produced in \cite{kuznetsov_homological_2007} make up a left Lefschetz decomposition.
Secondly, we call the category $\cD^{\vee}$ a HP dual, whereas in \cite{kuznetsov_homological_2007}, that term is reserved for a variety $Y$ (which in general then neither exists nor is unique) such that $D^{b}(Y) \cong \cD^{\vee}$.

Our proof allows us to drop some hypotheses in Thm.~\ref{thm:fundamentalTheoremOfHPD} as compared to \cite{kuznetsov_homological_2007}.
Firstly we drop the assumption that $\cL$ is globally generated, so that the starting data for HP duality is not a map $X \to \PP(H^{0}(X,\cL)^{\vee})$, but rather a pair of $(X,\cL)$, where $\cL$ is an effective line bundle.\footnote{Another approach to handling $\cL$ which are not globally generated is to pass to a blow-up of $X$, see \cite{carocci_homological_2015}.}
In the original set-up for HP duality, Kuznetsov shows that  $(\cD^{\vee})^{\vee} \cong \cD$, so that HPD is indeed a duality \cite[Thm.~7.3]{kuznetsov_homological_2007}.
We don't show this result here, but remark that it requires stronger hypotheses than we're taking; in particular $n \le \dim V$ and $\cC_{V} = 0$ (which follows from $\cL$ globally generated) seems necessary.

We also drop the assumption that the base changed varieties have the expected dimension by interpreting all base changes in a suitable derived sense in general.
Admittedly, from the point of view of geometry these non-proper derived intersections seem hard to understand.
Finally, instead of smooth projective varieties $X$ and $Y$, we work with saturated subcategories of smooth, quasi-projective quotient stacks.

This last generalisation is needed in this paper's main application, which is to the constructions of \cite{rennemo_homological_2015,rennemo_hori-mological_2016}.
In these papers, Ed Segal and I produce candidate HP duals for $\Sym^{2}\PP^{n}$ and the generalised Pfaffian varieties $\Pf(k,2n+1)$.
These candidate HP duals are obtained via similar VGIT for LG model techniques, which lets us prove roughly half of Thm.~\ref{thm:fundamentalTheoremOfHPD} directly.
Upgrading this to the full statement of Thm.~\ref{thm:fundamentalTheoremOfHPD} is tricky.
The result of this paper resolves this issue and lets us prove that the candidate HP duals are in fact HP duals.

\subsection{Acknowledgements}
I thank D.\,Beraldo, A.\,Kuznetsov, E.\,Segal and R.\,Thomas for helpful discussions related to this paper.

\subsection{Conventions}
Since semi-orthogonality conditions are insensitive to cohomological degree, we generally omit cohomological shifts in formulas.
All functors are derived.
The space $\Hom(\cE,\cF)$ is cohomologically graded and contains all shifted maps, i.e.~if $\cE$ and $\cF$ are sheaves, then $\Hom(\cE,\cF) = \bigoplus_{i} \Ext^{i}(\cE,\cF)$.

We work over $\CC$.

\section{Technical background}
\subsection{Saturated and admissible categories}
Let us recall some definitions and results about semi-orthogonal decompositions, all taken from \cite[Sec.~2]{kuznetsov_homological_2007}.
A semi-orthogonal decomposition of a triangulated category $\cC$ is a sequence $\cC_{1}, \ldots, \cC_{n}$ of full, triangulated subcategories of $\cC$ such that $\cC_{i} \subseteq \cC_{j}^{\perp}$ if $j > i$, and such that
\[
\cC = \langle \cC_{1}, \ldots, \cC_{n}\rangle,
\]
i.e.~the categories $\cC_{i}$ generate $\cC$.

A full, triangulated subcategory $\cC \subseteq \cD$ is called left (resp.~right) admissible if the inclusion functor admits a left (resp.~right) adjoint, and is called admissible if it is both left and right admissible.
\begin{nlemma}
  If $\cD \subseteq \cC$ is left admissible, then $\cC = \langle \cD, ^{\perp}\!\!\cD \rangle$ is a semi-orthogonal decomposition.
If $\cD \subseteq \cC$ is right admissible, then $\cC = \langle \cD^{\perp}, \cD \rangle$ is a semi-orthogonal decomposition. 
\end{nlemma}

A triangulated category $\cC$ is saturated if every exact functor $\cC \to D^{b}(\CC)$ and $\cC^{\op} \to D^{b}(\CC)$ is representable.
By a result of Bondal and Van den Bergh, the category $D^{b}(X)$ is saturated for a smooth, projective variety $X$ \cite{bondal_generators_2003}.
\begin{nlemma}
\label{thm:saturatedMeansAdmissible}
If $F \colon \cC \to \cD$ is a fully faithful inclusion of triangulated categories, $\cC$ is saturated and $\cD$ is ext finite, then $F$ admits both adjoints, i.e.~$\cC$ is an admissible subcategory.
\end{nlemma}

\begin{nlemma}
If $\cC \subseteq \cD$ is a left or right admissible subcategory and $\cD$ is saturated, then $\cD$ is saturated.
\end{nlemma}

\begin{nlemma}
\label{thm:compositionOfAdmissibles}
  If $\cC_{1}, \ldots, \cC_{n} \subseteq \cC$ is a sequence of admissible subcategories such that $\cC_{i} \subseteq \cC_{j}^{\perp}$ if $j > i$, then
\[
\cC = \langle \langle \cC_{1}, \ldots, \cC_{n} \rangle^{\perp}, \cC_{1}, \ldots, \cC_{n} \rangle
\]
and 
\[
\cC = \langle \cC_{1}, \ldots, \cC_{n}, ^{\perp}\!\!\langle \cC_{1}, \ldots, \cC_{n} \rangle \rangle
\]
are semi-orthogonal decompositions.
\end{nlemma}

\subsection{Matrix factorisation categories}
\label{sec:MFCategories}
We give a brief introduction to these categories to fix notation, and refer to \cite{addington_pfaffian-grassmannian_2015, ballard_variation_2012, rennemo_homological_2015, shipman_geometric_2012} for further background.
A Landau--Ginzburg (LG) model is for us the following data:
\begin{itemize}
\item A quotient stack $X$ of the form $[A/G \times \CC^{*}_{R}]$, where $G$ is a reductive algebraic group and $A$ is a smooth, quasi-projective variety.
\item A function $W \colon [A/G] \to \CC$, which is of degree 2 with respect to the $\CC^{*}_{R}$-action.
\end{itemize}
We require that $-1 \in \CC^{*}_{R}$ acts trivially on $[A/G]$, i.e.~that there exists a $g \in G$ such that $(g,-1) \in G \times \CC^{*}_{R}$ acts as the identity on $A$.

From this data, work of Positselski and Orlov shows that we can define a category of matrix factorisations $D^{b}(X,W)$ \cite{efimov_coherent_2015,orlov_matrix_2012, positselski_two_2011}.
Let $\cO_{X}[1]$ denote the line bundle corresponding to the fundamental character of $\CC^{*}_{R}$, and let $\cE[i] = \cE \otimes \cO_{X}[1]^{\otimes i}$ for any sheaf $\cE$ on $X$ -- the operation $\cE \mapsto \cE[1]$ will be the shift functor in the category $D^{b}(X,W)$.
The objects of the category $D^{b}(X,W)$ are pairs $(\cE,d)$, where $\cE$ is a coherent sheaf on $X$ and $d \colon \cE \to \cE[1]$ is a homomorphism such that $d^{2} = \id_{\cE} \otimes W$.
Given two objects $(\cE,d)$ and $(\cE',d')$, the space $\Hom(\cE,\cE') = \oplus_{i\in \ZZ} \Hom(\cE,\cE'[i])$ then obtains a differential, and so the category of such pairs $(\cE,d)$ becomes a dg category.
Taking the Verdier quotient of its homotopy category with respect to a subcategory of ``acyclic factorisations'', we obtain the category $D^{b}(X,W)$.

One computationally simple way to describe the hom spaces in this category is to note that if $(\cE, d)$ and $(\cE', d')$ are objects, then $\sHom(\cE,\cE')$ is a complex of sheaves on $X$, and if $\cE$ is locally free, then $\Hom_{D^{b}(X,W)}(\cE,\cE') \cong \RGamma(X, \sHom(\cE,\cE'))$.

It turns out that the category $D^{b}(X,W)$ is a natural generalisation of the usual derived categories:
If we let the $\CC^{*}_{R}$-action be trivial and set $W = 0$, we get
\begin{equation}
\label{eqn:MFVsDerived}
D^{b}(X,0) \cong D^{b}([A/G]),
\end{equation}
see \cite[Prop.~2.1.6]{ballard_homological_2013}.

We have the usual derived functors between factorisation categories, e.g.~if $f \colon (X,W) \to (X',W')$ is a morphism of LG models (meaning a map commuting with the $\CC^{*}_{R}$-actions and the potential), then one has a pullback map $f^{*}$ between the derived categories, and a push-forward map $f_{*}$ if $f$ is a closed immersion.

\begin{nremark}
Given the starting data of $A$ acted on by $G \times \CC^{*}_{R}$ and a potential $W$, we are free to modify the $\CC^{*}_{R}$-action in the following way.
Take a homomorphism $\phi \colon \CC^{*}_{R} \to G$, let $\sigma_{R} \colon \CC^{*}_{R} \to \Aut(A)$ and $\sigma_{G} \colon G \to \Aut(A)$ be the original actions, and define new actions by
\[
(\sigma_{G}', \sigma_{R}') = (\sigma_{G}, \sigma_{R} \cdot \sigma_{G} \circ \phi). 
\]
This does not change the category of matrix factorisations.
In particular, every time a category with potential $W = 0$ appears in the proof, we can apply this operation and assume it has the trivial $\CC^{*}_{R}$-action.
For this reason, and to align with (\ref{eqn:MFVsDerived}), we'll from this point on drop the $\CC^{*}_{R}$-action from the notation and call $(X,W) = ([A/G],W)$ an LG model.
\end{nremark}

\subsubsection{Kn\"orrer periodicity}
For a general potential $W$, the matrix factorisations categories are hard to analyse geometrically.
There is, however, one case in which they are well understood, via the result known as (global) Kn\"orrer periodicity.

Let $X$ be a smooth, quasi-projective quotient stack, and let $E \to X$ be a vector bundle.
Choose a section $s \colon X \to E$, and let $Y \subset X$ be the vanishing locus of $s$.
The section $s$ induces a function $W = s^{\vee} \colon E^{\vee} \to \CC$.
If we let $\CC^{*}_{R}$ act on $E^{\vee}$ by scaling the fibres with weight 2, we get an LG model $(E^{\vee},W)$.
Let $\pi \colon E^{\vee} \to \cX$ be the projection, let $i \colon \pi^{-1}(Y) \into E^{\vee}$ be the inclusion, and write $\pi$ also for the restricted projection $\pi^{-1}(Y) \to Y$.
\begin{nprop}
\label{thm:KnorrerPeriodicity}\cite{hirano_derived_2016,isik_equivalence_2013,shipman_geometric_2012}
If $Y$ has the expected codimension (i.e.~rank $E$) in $X$, then we have an equivalence
\[
\pi_{*}\circ i^{*} \colon D^{b}(E^{\vee},W) \to D^{b}(Y).
\]
\end{nprop}
So any scheme defined by the vanishing of a regular section of a vector bundle is derived equivalent to a factorisation category in a natural way.
If the dimension of $Y$ is higher than expected, then the category $D^{b}(E^{\vee},W)$ is likely equivalent to the derived category of coherent sheaves on the derived zero locus of $s$; see e.g.~\cite[8.2.2]{arinkin_singular_2015} for a definition of this category.
At any rate, the results in this paper indicate that $D^{b}(E^{\vee},W)$ is the correct replacement of $D^{b}(Y)$ for the purposes of HP duality.

\subsection{VGIT and matrix factorisation categories}
Our arguments rely on results about variation of GIT stability for matrix factorisation categories, going back to Segal for a $\CC^{*}$-quotient (which is the case we need) \cite{segal_equivalence_2011}, and worked out by Halpern-Leistner and Ballard--Favero--Katzarkov for a general GIT quotient \cite{halpern-leistner_derived_2015,ballard_variation_2012}.

We'll only need a very special case of the general theory, concerning a special kind of $\CC^{*}$-action.
Let $(X,W) = [A/G]$ be a smooth quotient stack, let $E \to X$ be a vector bundle, and let $\CC^{*}$ act linearly on $E$, fixing $X$.
Let $Z$ be the total space of $E$, and choose a $\CC^{*}_{R}$-action and a $\CC^{*}$-invariant potential $W$ on $Z$.

Choosing the positive or negative character of $\CC^{*}$ gives two GIT stability conditions for the quotient problem of $\CC^{*}$ acting on $Z$.
The unstable loci are subbundles $Y_{\pm} \subseteq Z$, where $Y_{+}$ (resp.~$Y_{-}$) is the sub-bundle on which $\CC^{*}$ acts with non-positive (resp.~non-negative) weights.
The stable loci are the complements $Z_{\pm} = Z \setminus Y_{\pm}$.

Let $N_{\pm}$ be the normal bundle of $Y_{\pm}$ in $Z$.
The line bundle $\wedge^{\dim N_{\pm}} N_{\pm}|_{X}$ carries a $\CC^{*}$-action, and since $X$ is $\CC^{*}$-fixed, the weight of the $\CC^{*}$-action on this line bundle is well defined.
Let $n_{\pm}$ denote $\pm 1$ times this weight.

Now, again since $X$ is $\CC^{*}$-fixed, any object $\cE \in D^{b}([X/\CC^{*}],W)$ splits into eigensheaves for the $\CC^{*}$-action, and this gives an orthogonal decomposition
\[
D^{b}([X/\CC^{*}]) = \oplus_{i\in \ZZ} D^{b}(X)(i),
\]
where $(i)$ denotes twisting by $i$ times the identity character of $\CC^{*}$.
Define $\cW_{[0,n]} \subseteq D^{b}([Z/\CC^{*}],W)$ to be the full subcategory consisting of those objects $\cE$ such that the restriction of $\cE$ to $D^{b}([X/\CC^{*}])$ lies in 
\[
\bigoplus_{0 \le i \le n} D^{b}(X)(i).
\]
The main result of \cite{ballard_variation_2012,halpern-leistner_derived_2015} then specialises to the following claim.
\begin{nprop}
\label{thm:GITResult}
The restriction functor $D^{b}([Z/\CC^{*}],W) \to D^{b}([Z_{\pm}/\CC^{*}],W)$ induces an equivalence $\cW_{[0,n_{\pm} -1]} \to D^{b}([Z_{\pm}/\CC^{*}],W)$.
\end{nprop}

\section{Proofs}
\label{sec:HPDSection}
Let $X$ be a stack of the form $[A/G]$, where $A$ is a smooth, quasiprojective variety and $G$ is a reductive algebraic group acting on $A$.
We assume that $X$ is ext finite, i.e.~that for any two coherent sheaves $\cE, \cF$ on $X$ the space $\oplus_{i \in \ZZ} \Ext^{i}(\cE,\cF)$ is finite-dimensional.
Let $\cL$ be a line bundle on $X$, and let $V = H^{0}(X,\cL)^{\vee}$.

Let $\cD \subseteq D^{b}(X)$ be a full, triangulated subcategory.
We assume that $\cD$ is closed under the operation $\cE \mapsto \cE \otimes \cL^{\pm 1}$, and that we are given a Lefschetz decomposition
\[
\cD = \langle \cA_{0}, \cA_{1}(\cL), \ldots, \cA_{n}(n\cL) \rangle,
\]
with $\cA_{0} \supseteq \cA_{1} \supseteq \cdots \supseteq \cA_{n}$.
We furthermore assume that all the $\cA_{i}$ are saturated categories.

\subsection{The GIT quotients}
Fix a linear subspace $L \subseteq V^{\vee}$ of dimension $l$.
Let $\wt{X}$ be the total space of $\cL^{\vee}$, and let $Z = \wt{X} \times L$.
Equip $Z$ with the $\CC^{*}$-action that scales the fibres of $\wt X$ with weight 1 and scales $L$ with weight $-1$, and let $\cZ = [\wt X \times L/\CC^{*}]$.
To be precise about our conventions, this means that the line bundle $\cO_{\cZ}(1)$ has no sections when restricted to $L$, while $\cO_{\cZ}(-1)$ has no sections when restricted to $\wt X$.

We take the $\CC^{*}_{R}$-action on $\cZ$ which scales the fibres of $\wt X \to X$ with weight 2 and fixes $L$, and the potential
\[
W \colon \cZ = \wt X \times_{\CC^{*}} L \into \wt X \times_{\CC^{*}} H^{0}(\cX,\cL) \to \CC,
\] 
where the last map is the obvious pairing. 
This makes $(\cZ, W)$ into an LG model.

There are two natural GIT stability conditions for the action of $\CC^{*}$ on $Z$.
For the first, ``positive'' stability, the unstable locus is $Y_{+} = X \times L \subset Z$.
For the second, ``negative'' stability, the unstable locus is $Y_{-} = \wt{X} \times \{0\} \subset Z$.
We denote the semistable loci by $Z_{\pm} = Z \setminus Y_{\pm}$.

For future reference, we let $i_{\pm} \colon Y_{\pm} \to Z$, $i_{0} \colon X \into Z$ and $j_{\pm} \colon Z_{\pm} \to Z$ denote the inclusions, and we let $\pi_{\pm} \colon Y_{\pm} \to X$ denote the projections.
For $X, \wt X, Y_{\pm}, Z_{\pm}$, we use the calligraphic versions of the same letters to denote the stacks obtained by quotienting out by the $\CC^{*}$-action, e.g.~$\cX = [X/\CC^{*}]$.

The $\CC^{*}$-action on $X$ is trivial, so we have a decomposition
\[
D^{b}(\cX) \cong \bigoplus_{i \in \ZZ} D^{b}(X)(i).
\]
Given any subcategory $\cC \subseteq D^{b}(X)$, we'll write $\cC(i)$ for the corresponding subcategory of $D^{b}(X)(i) \subseteq D^{b}(\cX)$.

The stack $\cZ$ has two natural line bundles on it: The line bundle $\cL$ pulled back from $\cX$, and the line bundle $\cO(1)$, corresponding to the identity character of $\CC^{*}$.
If $\cE$ is a sheaf on $\cZ$, we write $\cE(i\cL,j)$ for $\cE \otimes \cL^{\otimes i} \otimes \cO(1)^{\otimes j}$.

\subsection{Defining the categories $\cD_{L^{\perp}}$ and $\cD_{L}^{\vee}$}
\label{sec:baseLocusCategory}
Let $X_{L^{\perp}} \subseteq X$ be the base locus of the linear system $\PP L^{\perp} \subset \PP V^{\vee} = |\cL|$, i.e.~the substack cut out by the map 
\[
L \otimes \cL^{\vee} \into H^{0}(X,\cL) \otimes \cL^{\vee} \to \cO_{X}.
\]
We let $\cD_{L^{\perp}} \subseteq D^{b}(\cZ_{+},W)$ be the subcategory of objects $\cE$ such that $i_{0}^{*}\cE \in \cD$.
The category $\cD_{L^{\perp}}$ should be thought of as the (derived) restriction of $\cD$ to $X_{L^{\perp}}$.

Note that $\cZ_{+}$ is naturally isomorphic to the total space of the vector bundle $\cL^{\vee} \otimes L$ on $X$, so that if $X_{L^{\perp}}$ has the expected dimension, we have the Kn\"orrer equivalence of Prop.~\ref{thm:KnorrerPeriodicity}:
\begin{equation}
\label{eqn:KnorrerEquivalence}
\pi_{*}\circ i^{*} \colon D^{b}(X_{L^{\perp}}) \stackrel{\cong}{\to} D^{b}(\cZ_{+},W).
\end{equation}
The definition of $\cD_{L^{\perp}}$ is motivated by the following easy proposition.
\begin{nprop}
If $X_{L^{\perp}}$ has the expected dimension, then under the equivalence \eqref{eqn:KnorrerEquivalence}, the category $\cD_{L^{\perp}}$ becomes the category of those $\cE \in D^{b}(X_{L^{\perp}})$ such that the push-forward of $\cE$ along $D^{b}(X_{L^{\perp}}) \to D^{b}(X)$ lies in $\cD$.
\end{nprop}

The category $\cD_{L}^{\vee}$ is a full subcategory of $D^{b}(\cZ_{-},W)$, defined as follows.
The stack $\cZ_{-}$ has a natural map to $\PP L$, and the fibre at the point $[f] \in \PP L \subseteq \PP(H^{0}(X,\cL))$ is isomorphic to $(\wt X, f)$.
For any $[f]$, we have an inclusion 
\[
X = X \times [f] \into \wt X \times_{\CC^{*}} (L \setminus 0) \into \cZ_{-}.
\]
We let $\cD_{L}^{\vee}$ be the subcategory of $D^{b}(\cZ_{-},W)$ consisting of those objects $\cE$ such that for each $[f] \in \PP L$ we have 
\[
\cE|_{X \times [f]} \in \cA_{0}.
\]

\subsection{The windows}
\label{sec:windows}
Given any full, triangulated subcategory $\cS \subseteq D^{b}(\cX)$, we define the ``window subcategory'' $\cW(\cS) \subseteq D^{b}(\cZ,W)$ as the subcategory of objects $\cE$ such that $i_{0}^{*}\cE \in \cS$.

We define two subcategories $\cS_{\pm} \subseteq D^{b}(\cX)$ by $\cS_{+} = \langle \cA_{i}(i) \rangle_{0 \le i \le n}$, and $\cS_{-} = \langle \cA_{0}(i) \rangle_{0 \le i \le l-1}$.
The relation between $\cD_{L^{\perp}}$ and $\cD^{\vee}_{L}$ is obtained from the following proposition, which is proved by combining Lemmas \ref{thm:fullyFaithfulConditionPlus} and \ref{thm:windowEquivalencePlusCase} with Cor.~\ref{thm:windowEquivalenceMinus}.
\begin{nprop}
\label{thm:OurSpecialGITStatement}
The restriction functors induce equivalences 
\[
j_{+}^{*} \colon \cW(\cS_{+}) \to \cD_{L^{\perp}}
\]
and
\[
j_{-}^{*} \colon \cW(\cS_{-}) \to \cD_{L}^{\vee}
\]
\end{nprop}

We further let $\cS_{0} = \cS_{+} \cap \cS_{-} = \langle \cA_{i}(i) \rangle_{0 \le i \le \min(l-1,n)}$ and define $\cC_{L} = \cW(\cS_{0})$.
Given Prop.~\ref{thm:OurSpecialGITStatement}, it follows that $\cC_{L}$ embeds into both $\cD_{L^{\perp}}$ and $\cD^{\vee}_{L}$.

Let's record for future use the following specialisation of Prop.~\ref{thm:GITResult}.
\begin{nprop}
\label{thm:GITResultOurCase}
The restriction functors induce equivalences 
\[
j_{+}^{*} \colon \cW(D^{b}(X)(0)) \to D^{b}(\cZ_{+},W)
\]
and
\[
j_{-}^{*} \colon \cW(\langle D^{b}(X)(0), \ldots, D^{b}(X)(l-1) \rangle) \to D^{b}(\cZ_{-},W).
\]
\end{nprop}
\begin{remark}
It is a key point in our proof that the subcategory $\cS_{+} \subseteq D^{b}(\cX)$ used to define $\cW(\cS_{+})$ is not the same as the one produced by the general theory of \cite{ballard_variation_2012, halpern-leistner_derived_2015}.
Following the lines of the general theory we could consider 
\begin{equation}
\label{eqn:standardWindow}
\cS_{+}^{\std} = \cD(0) = \langle \cA_{0}(0), \cA_{1}(0,\cL), \ldots, \cA_{n}(0,n\cL)\rangle \subset D^{b}(\cX).
\end{equation}
It's a simple corollary of the general theory (as summarised in Prop.~\ref{thm:GITResultOurCase}) that $\cW(\cS_{+}^{\std}) \cong \cD_{L^{\perp}}$, but there is no obvious way to compare $\cW(\cS_{+}^{\std})$ with $\cW(\cS_{-})$.

Instead our $\cS_{+} = \langle \cA_{0}, \cA_{1}(1), \ldots, \cA_{n}(n)\rangle$ is obtained by twisting the pieces of the decomposition (\ref{eqn:standardWindow}), so that $\cW(\cS_{+})$ realises $\cD_{L^{\perp}}$ as a subcategory of $D^{b}(\cZ,W)$ in a non-standard way.
This means that there is extra work to be done in proving $\cW(\cS_{+}) \cong \cD_{L^{\perp}}$, but the reward for this is that $\cW(\cS_{+}) \cap \cW(\cS_{-})$ is (usually) non-empty, giving the relation claimed in Thm.~\ref{thm:fundamentalTheoremOfHPD} between $\cD_{L^{\perp}}$ and $\cD_{L}^{\vee}$.
\end{remark}

\subsubsection{Fully faithfulness of $\cW(\cS_{+}) \to \cD_{L^{\perp}}$}
Let us now give a few results about when $\Hom(\cE,\cF) \to \Hom(j_{\pm}^{*}\cE, j_{\pm}^{*}\cF)$ is an equivalence.
\begin{nlemma}
\label{thm:localCohomology}
For $\cE, \cF \in D^{b}(\cZ,W)$, the restriction map $\Hom(\cE, \cF) \to \Hom(j_{+}^{*}\cE, j_{+}^{*}\cF)$ is an isomorphism if $\Hom(i_{+}^{*}\cE, i_{+}^{*}\cF(-i\cL,i)) = 0$ for all $i > 0$.

For $\cE, \cF \in D^{b}(\cZ,W)$, the restriction map $\Hom(\cE, \cF) \to \Hom(j_{-}^{*}\cE, j_{-}^{*}\cF)$ is an isomorphism if $\Hom(i_{-}^{*}\cE, i_{-}^{*}\cF(i)) = 0$ for all $i \le -l$.
\end{nlemma}
\begin{proof}
Let $N_{\pm}$ be the normal bundle to $\cY_{\pm}$ in $\cZ$. 
The local cohomology argument of the proof of \cite[Lemma 5.14]{rennemo_homological_2015} (which appears in many places, e.g.~\cite{ballard_variation_2012, halpern-leistner_derived_2015, segal_equivalence_2011, teleman_quantization_2000}) shows that the restriction map $\Hom_{\cZ}(\cE,\cF) \to \Hom_{\cZ_{\pm}}(\cE|_{\cZ_{\pm}}, \cF|_{\cZ_{\pm}})$ is fully faithful if $\Hom(\cE|_{\cY_{\pm}},\cF|_{\cY_{\pm}} \otimes \wedge^{\dim N_{\pm}} N_{\pm} \otimes \Sym^{i}N_{\pm}) = 0$ for $i \ge 0$.
The claims now follow from $N_{+} = \cL^{\vee}(1)$ and $N_{-} \cong \cO(-1)^{\oplus l}$.
\end{proof}

The following lemma is standard.
\begin{nlemma}
\label{thm:SODOfL}
\label{thm:SODOfTildeX}
The functors $\pi_{\pm}^{*} \colon D^{b}(X)(i) \to D^{b}(\cY_{\pm})$ are fully faithful for all $i \in \ZZ$, and there are infinite semi-orthogonal decompositions
\[
D^{b}(\cY_{-}) = \langle \ldots,\pi_{-}^{*}D^{b}(X)(i-1), \pi_{-}^{*}D^{b}(X)(i), \pi_{-}^{*}D^{b}(X)(i+1), \ldots \rangle
\]
and
\[
D^{b}(\cY_{+}) = \langle \ldots,\pi_{+}^{*}D^{b}(X)(i+1), \pi_{+}^{*}D^{b}(X)(i), \pi_{+}^{*}D^{b}(X)(i-1), \ldots \rangle.
\]
\end{nlemma}

\begin{nlemma}
\label{thm:HelpfulTrick}
If $\cE, \cF \in D^{b}(\cY_{+})$ are such that $\Hom(i_{0}^{*}\cE,i_{0}^{*}\cF(i)) = 0$ for $i \ge 1$, then the restriction map $\Hom(\cE, \cF) \to \Hom(i_{0}^{*}\cE, i_{0}^{*}\cF)$ is an isomorphism.

If $\cE, \cF \in D^{b}(\cY_{-})$ are such that $\Hom(i_{0}^{*}\cE,i_{0}^{*}\cF(i)) = 0$ for $i \le -1$, then the restriction map $\Hom(\cE, \cF) \to \Hom(i_{0}^{*}\cE, i_{0}^{*}\cF)$ is an isomorphism.
\end{nlemma}
\begin{proof}
We only treat the ``+'' case, the ``-'' case is the same.
We may write $\cE$ as an iterated extension of objects $\pi_{+}^{*}\cE_{i}(i)$, where $\cE_{i} \in D^{b}(X)(0)$, and similarly for $\cF$.
The condition $\Hom(i_{0}^{*}\cE,i_{0}^{*}\cF(i)) = 0$ for $i \ge 1$ implies that $\Hom(\cE_{i}, \cF_{j}) = 0$ for $j < i$.
Thus $\Hom(\pi_{+}^{*}\cE_{i}(i), \pi_{+}^{*}\cF_{j}(j)) = \Hom(\cE_{i}, \cF_{j} \otimes (\pi_{+})_{*}(\cO_{\cY_{+}}(j-i))) = 0$.
Hence $\Hom(\pi_{+}^{*}\cE_{i}(i), \pi_{+}^{*}\cF_{j}(j)) = 0$ unless $i = j$.

Let $i \in \ZZ$ be minimal such that $\cE_{i}$ or $\cF_{i}$ is non-trivial.
We have exact triangles $\pi_{+}^{*}\cE_{i}(i) \to \cE \to \cE_{< i}$ and $\pi_{+}^{*}\cF_{i}(i) \to \cF \to \cF_{< i}$.
Now $\Hom(\pi_{+}^{*}\cE_{i}(i), \cF_{<i}) = 0$ and $\Hom(\cE_{<i}, \pi_{+}^{*}\cF_{i}(i)) = 0$, and drawing up all the associated exact triangles of hom spaces, this implies that $\Hom(\cE,\cF) = \Hom(\pi^{*}_{+}\cE_{i}(i), \pi_{+}^{*}\cF_{i}(i)) \oplus \Hom(\cE_{<i}, \cF_{<i})$.
Hence by induction $\Hom(\cE,\cF) = \oplus_{i \in \ZZ} \Hom(\pi^{*}_{+}\cE_{i}(i), \pi_{+}^{*}\cF_{i}(i)) = \Hom(i_{0}^{*}\cE,i_{0}^{*}\cF)$.
\end{proof}

\begin{nlemma}
\label{thm:bestLocalCohomologyStatement}
For $\cE, \cF \in D^{b}(\cZ,W)$, the restriction map
\[
\Hom(\cE, \cF) \to \Hom(j_{+}^{*}\cE, j_{+}^{*}\cF)
\]
is an isomorphism if
\[
\Hom(i_{0}^{*}\cE, i_{0}^{*}\cF(-i\cL, i + j)) = 0
\]
for all $i > 0, j \ge 0$.

For $\cE, \cF \in D^{b}(\cZ,W)$, the restriction map
\[
\Hom(\cE, \cF) \to \Hom(j_{-}^{*}\cE, j_{-}^{*}\cF)
\]
is an isomorphism if $\Hom(i_{0}^{*}\cE, i_{0}^{*}\cF(i)) = 0$ for all $i \le -l$.
\end{nlemma}
\begin{proof}
Combine Lemmas \ref{thm:localCohomology} and \ref{thm:HelpfulTrick}.
\end{proof}

\begin{nlemma}
\label{thm:fullyFaithfulConditionPlus}
If $\cE, \cF \in \cW(\cS_{+})$, then $\Hom(\cE, \cF) \to \Hom(\cE|_{\cZ_{+}}, \cF|_{\cZ_{+}})$ is an isomorphism.
\end{nlemma}
\begin{proof}
By Lemma \ref{thm:bestLocalCohomologyStatement}, it's enough to check that if $\cE, \cF \in \cW(\cS_{+})$ we have 
\[
\Hom(i_{0}^{*}\cE,i_{0}^{*}\cF(-i\cL,i+j)) = 0
\]
for all $i > 1$ and $j \ge 0$.
Since $i_{0}^{*}(\cE)$ and $i_{0}^{*}(\cF)$ lie in $\cS_{+} = \langle \cA_{i}(i) \rangle_{i \in [0,n]}$, it suffices to check this under the assumption that $i_{0}^{*}\cE \in \cA_{k}(k)$ and $i_{0}^{*}\cF \in \cA_{k'}(k')$ for some $k, k' \in [0,n]$.
The claim is then a simple consequence of the Lefschetz semi-orthogonality property for the $\cA_{i}$.
\end{proof}

\subsubsection{Essential surjectivity of $\cW(\cS_{+}) \to \cD_{L^{\perp}}$}
A simple amplification of Prop.~\ref{thm:GITResultOurCase} shows that $j_{+}^{*}\colon \cW(\cD(0)) \to \cD_{L^{\perp}}$ is an equivalence.
So it suffices to show that every object of $\cW(\cS_{+})$ is equivalent to an object of $\cW(\cD(0))$, up to taking cones over objects supported on $\cY_{+}$.
This follows from the more precise result that $\cW(\langle \cD(i) \rangle_{i \in \ZZ}) \subseteq D^{b}(\cZ,W)$ admits a certain semiorthogonal decomposition, see Lemma \ref{thm:SODOfZ}.

It will be useful in calculations to use the notation $(i_{+})_{!} = (i_{+})_{*} \otimes \cL^{\vee}(1)$, which is justified by the fact that $(i_{+})_{!}$ is left adjoint to $i_{+}^{*}$.
Define the functor $\Phi = (i_{+})_{!}\pi_{+}^{*} \colon D^{b}(\cX) \to D^{b}(\cZ,W)$.
\begin{nlemma}
\label{thm:InclusionIntoZMinusFullyFaithful}
The functor $\Phi \colon D^{b}(X)(k) \to D^{b}(\cZ,W)$ is fully faithful for any $k \in \ZZ$.
\end{nlemma}
\begin{proof}
We may WLOG assume $k = 0$, so let $\cE, \cF \in D^{b}(X)(0)$. We have $\Hom(\Phi\cE, \Phi\cF) = \Hom(\pi_{+}^{*}\cE, i_{+}^{*}(i_{+})_{!}\pi_{+}^{*}\cF)$, and the counit $i_{+}^{*}(i_{+})_{!} \to \id_{[L/\CC^{*}]}$ induces an exact triangle
\[
\pi_{+}^{*}\cF(-\cL,1) \to i_{+}^{*}(i_{+})_{!}\pi_{+}^{*}\cF \to \pi_{+}^{*}\cF.
\]
Now
\[
\Hom(i_{0}^{*}\pi_{+}^{*}\cE, i_{0}^{*}\pi_{+}^{*}(\cF \otimes \cL^{\vee}(1)) \otimes \cO(i)) = \Hom(\cE, \cF(-\cL,i+1)) = 0
\] 
for all $i \ge 0$.
It follows by Lemma \ref{thm:HelpfulTrick} that $\Hom(\pi_{+}^{*}\cE, \pi_{+}^{*}\cF(-\cL,1)) = 0$, and so $\Hom(\cE,\cF) = \Hom(\Phi\cE, \Phi\cF)$ as required.
\end{proof}

\begin{nlemma}
\label{thm:SODOfZ}
There exists a semi-orthogonal decomposition
\[
\cW(\langle \cD(k)\rangle_{k \in \ZZ}) \cong \langle \cW(\langle \cA_{0}(k) \rangle_{k \in \ZZ}), \langle \Phi\cA_{i}(i\cL,j) \rangle_{i,j} \rangle,
\]
where $i = 1, \ldots, n$ and $j \in \ZZ$.
The semi-orthogonality relations are such that $\Phi\cA_{i}(i\cL,j) \subseteq \Phi\cA_{i'}(i'\cL,j')^{\perp}$ if
\begin{itemize}
\item $i < i'$, or if
\item $i = i'$ and $j > j'$.
\end{itemize}
\end{nlemma}
\begin{proof}
The restriction of the functor $\Phi$ to $D^{b}(X)(0)$ admits a right adjoint given by $(-)^{\CC^{*}}\circ(\pi_{+})_{*}(i_{+})^{*}$.
Hence the restriction of $\Phi$ to $\cA_{i}(i\cL,j)$ admits a right adjoint for all $i$ and $j$, and so the categories $\Phi\cA_{i}(i\cL,j)$ are right admissible.

It's therefore enough to prove the semi-orthogonality properties and
\begin{equation}
\label{eqn:cWistheSOD}
\cW(\langle \cA_{0}(k) \rangle_{k \in \ZZ}) = \langle \Phi\cA_{i}(i\cL,j) \rangle_{(i,j) \in [1,n] \times \ZZ}^{\perp}.
\end{equation}

Let $\cE \in \cA_{i}(i\cL,j)$ and $\cF \in \cA_{i'}(i'\cL,j')$ with $i,i' \ge 1$.
If $i < i'$, we have
\[
\Hom(\cF, \cE(k)) = \Hom(\cF, \cE(-\cL,1 + k)) = 0
\]
for all $k \in \ZZ$ by the Lefschetz property.
Arguing as in the proof of Lemma \ref{thm:InclusionIntoZMinusFullyFaithful}, it follows that $\Hom(\Phi\cF, \Phi\cE) = 0$.
Similarly, if $i = i'$ and $j > j'$, we easily find
\[
\Hom(\cF, \cE(k)) = \Hom(\cF, \cE(-\cL,1 + k)) = 0
\]
for $k \ge 0$, and so again $\Hom(\Phi\cF, \Phi\cE) = 0$.
This proves the semi-orthogonality relations between the categories $\Phi\cA_{i}(i\cL,j)$.

It remains to prove (\ref{eqn:cWistheSOD}).
Let first $\cE \in \cW(\langle \cA_{0}(k) \rangle_{k \in \ZZ})$.
If $\cF \in \cA_{i}(i\cL,j)$ with $i \ge 1$, we use Lemma \ref{thm:HelpfulTrick} to find $\Hom(\Phi\cF,\cE) = \Hom(\pi_{+}^{*}\cF,i_{+}^{*}\cE) = 0$, since $\Hom(\cF, \cE|_{\cX}(k)) = 0$ for all $k \in \ZZ$.

For the converse, let $\cE \not\in \cW(\langle \cA_{0}(k) \rangle_{k \in \ZZ})$.
There is an obvious semiorthogonal decomposition
\[
D^{b}(\cX) = \bigoplus_{j \in \ZZ} D^{b}(X)(j) = \bigoplus_{j \in \ZZ} \langle \cA_{i}(i\cL,j) \rangle_{i \in [0,n]},
\]
Decomposing $\cE|_{\cX}$ according to this, let $(i,j)$ be the pair such that $\cE|_{\cX}$ has a non-trivial component $\cE_{i,j} \in \cA_{i}(i\cL,j)$, where we require that $i$ is maximal with this property, and that $j$ is then minimal with this property for the chosen $i$.
We then have $\Hom(\cE_{i,j}, \cE|_{\cX}) \not= 0$.

By our assumption on $\cE$, we have $i \ge 1$.
We have as before $\Hom(\Phi\cE_{i,j}, \cE) = \Hom(\pi_{+}^{*}\cE_{i,j}, i_{+}^{*}\cE)$.
Using the maximality of $i$ and minimality of $j$, we get $\Hom(\cE_{i,j}, \cE|_{\cX}(k)) = 0$ for $k > 0$, and since $\Hom(\cE_{i,j},\cE|_{\cX}) \not= 0$, it follows by Lemma \ref{thm:HelpfulTrick} that $\Hom(\Phi\cE_{i,j},\cE) \not= 0$.
But then $\cE \not\in \langle \Phi(\cA_{i}(i\cL,j)) \rangle^{\perp}_{(i,j) \in [1,n]\times \ZZ}$.
\end{proof}

Let $\iota \colon \cW(\langle\cA_{0}(i)\rangle_{i \in \ZZ}) \to \cW(\langle\cD(i)\rangle_{i \in \ZZ})$ be the inclusion, and let $\iota^{*}$ be the left adjoint, which exists by Lemma \ref{thm:SODOfZ}.
\begin{nlemma}
\label{thm:projectionOfWeights}
Let $\cE \in \cW(\langle \cA_{i}(i\cL,j) \rangle_{(i,j) \in I}$, where $I \subseteq [0,n] \times \ZZ$.
Then $\iota^{*}\cE \in \cW(\langle \cA_{i}(i + j)\rangle_{(i,j) \in I})$.

Let $I' = I \cap ([1,n] \times \ZZ)$, and let $J = \bigcup_{(i,j) \in I'} \{j, j+1, \ldots, j+i \}$.
The cone over $\iota^{*}\cE \to \cE$ lies in $\cW(\langle\cD(j)\rangle_{j \in J})$.
\end{nlemma}
\begin{proof}
The projected object $\iota^{*}\cE$ can be constructed as follows.
Let $(i,j) \in I$, let $i$ be maximal for this property, and let $j$ be minimal for the given $i$.
Let $\cE_{(i,j)}$ be the projection of $\cE|_{\cX}$ to $\cA_{i}(i\cL,j)$.
Then as in the proof of Lemma \ref{thm:SODOfZ} there is a map $\Phi\cE_{(i,j)} \to \cE$, with cone $C$.
Using the exact triangle
\[
\cE_{(i,j)} \to i_{0}^{*}\Phi\cE_{(i,j)} \to \cE_{(i,j)} \otimes \cL^{\vee}(1),
\]
one then checks that $C \in \cW(\langle \cA_{i}(i\cL,j) \rangle_{(i,j) \in J}$, where $J = (I \setminus (i,j)) \cup (i-1,j+1)$.
Replacing $\cE$ with $C$ and repeating this procedure, we eventually end up with $\cE \in \cW(\langle \cA_{i}(i + j)\rangle_{(i,j) \in I})$.
The claim about the cone of $\iota^{*}\cE \to \cE$ is easy to see.
\end{proof}

\begin{nlemma}
\label{thm:windowEquivalencePlusCase}
The restriction functor $j_{+}^{*} \colon \cW(\cS_{+}) \to D^{b}(\cZ_{+},W)$ has essential image $\cD_{L^{\perp}}$.
\end{nlemma}
\begin{proof}
By Prop.~\ref{thm:GITResultOurCase}, the functor $j_{+}^{*}$ induces an equivalence $\cW(D^{b}(X)(0)) \to D^{b}(\cZ_{+},W)$, and we first claim that this restricts to give $\cW(\cD(0)) \cong \cD_{L^{\perp}}$.
If $\cE \in \cW(D^{b}(X)(0))$, so that $i_{0}^{*}\cE \in D^{b}(X)(0)$, then by Lemma \ref{thm:SODOfTildeX} we must have $i_{-}^{*}\cE \in \pi_{-}^{*}D^{b}(X)(0)$.
It follows that $i_{-}^{*}\cE = \pi_{-}^{*}i_{0}^{*}\cE$.
Using this we find that the following diagram of functors commutes:
\[
\begin{tikzcd}
\cW(D^{b}(X)(0)) \arrow[r, "j_{+}^{*}"] \arrow[d, "i_{-}^{*}"] \arrow[dr, "i_{0}^{*}"] & D^{b}(\cZ_{+},W) \arrow[d,"|_{X}"] \\
D^{b}(\wt X) \arrow[r,"|_{\wt X \setminus X}"] & D^{b}(X)
\end{tikzcd}
\]
Thus $i_{0}^{*}\cE \in \cD(0)$ if and only if $(j_{+}^{*}\cE)|_{X} \in \cD$, proving our claim.

Hence any $\cE \in \cD_{L^{\perp}}$ is the restriction of some $\cE' \in \cW(\cD(0))$.
By Lemma \ref{thm:projectionOfWeights}, we have $\iota^{*}\cE' \in \cW(\cS_{+})$, and since $j_{+}^{*}\circ\Phi = 0$, we have $j_{+}\iota^{*}\cE' = j_{+}^{*}\cE' = \cE$.
Hence the functor $j_{+}^{*} : \cW(S_{+}) \to D^{b}(\cZ_{+},W)$ is essentially surjective, and Lemma \ref{thm:fullyFaithfulConditionPlus} shows that it is fully faithful.
\end{proof}

\subsubsection{Decomposition of $\cD_{L^\perp}$}
We now show that the semiorthogonal complement to $\cC_{L}$ inside $\cW(\cS_{+}) \cong \cD_{L^{\perp}}$ has the form claimed in Thm.~\ref{thm:fundamentalTheoremOfHPD}.

\begin{nlemma}
\label{thm:UsefulTrickVariant}
If $\cE, \cF \in D^{b}(\cX)$ and $\Hom(\cE, \cF(i\cL,-i)) = 0$ for $i \ge 1$, then the map $\Hom(\cE, \cF) \to \Hom(\pi_{-}^{*}\cE, \pi_{-}^{*}\cF)$ is an isomorphism.
\end{nlemma}
\begin{proof}
We have 
\[
\Hom(\pi_{-}^{*}\cE, \pi_{-}^{*}\cF) = \Hom(\cE, \cF \otimes (\pi_{-})_{*}\cO_{\wt \cX}) = \Hom(\cE, \cF \otimes \Sym^{\bullet}(\cL(-1))),
\]
and the claim follows.
\end{proof}
Let $(i_{-})_{!} = (i_{-})_{*} \otimes \cO(-l)$, so that $(i_{-})_{!}$ is left adjoint to $(i_{-})^{*}$, and define $\Psi = (i_{-})_{!}\pi_{-}^{*} \colon D^{b}(\cX) \to D^{b}(\cZ,W)$.
\begin{nlemma}
For $i \ge l$, the functor $\Psi \colon \cA_{i}(i) \to D^{b}(\cZ,W)$ is fully faithful and has image in $\cW(\cS_{+})$.
If $i > i' \ge l$, then $\Psi(\cA_{i}(i)) \in {}^{\perp}\Psi(\cA_{i'}(i'))$.
\end{nlemma}
\begin{proof}
Let $\pi \colon \cZ \to \cX$ be the projection.
For any $\cE \in D^{b}(\cX)$, we have $\Psi\cE = (i_{-})_{!}\pi_{-}^{*}\cE = \pi^{*}\cE \otimes (i_{-})_{!}\pi_{-}^{*}\cO_{\cX}$.
We have $(i_{-})_{!}\pi_{-}^{*}\cO_{\cX} = \cO_{\cY_{-}}(-l)$, and taking a deformed Koszul resolution of this sheaf (see e.g.~\cite[Sec.~3.3]{rennemo_homological_2015}) gives a locally free representation of $(i_{-})_{!}\pi_{-}^{*}\cO_{\cX}$ whose underlying sheaf is a direct sum of line bundles $\cO(-l), \ldots, \cO$. 
Thus $i_{0}^{*}(i_{-})_{!}\pi_{-}^{*}\cO_{\cX} \in \langle \cO_{\cX}(-l), \ldots, \cO_{\cX} \rangle$.
It follows that if $\cE \in \cA_{i}(i)$, then $i_{0}^{*}\Psi\cE \in \langle \cA_{i-l}(i-l), \ldots \cA_{i}(i)) \rangle \subseteq \cS_{+}$.

Let now $\cE \in \cA_{i}(i)$ and $\cE' \in \cA_{i'}(i')$, with $i \ge i' \ge l$.
It's enough to prove that $\Hom(\Psi\cE,\Psi\cE') = \Hom(\cE,\cE')$, since applying this with $i = i'$ gives fully faithfulness and applying it with $i > i'$ gives semi-orthogonality.

Adjunction gives
\[
\Hom(\Psi\cE, \Psi\cE') = \Hom(\pi_{-}^{*}\cE, i_{-}^{*}(i_{-})_{!}\pi_{-}^{*}\cE').
\]
In the exact triangle
\[
C \to i_{-}^{*}(i_{-})_{!}\pi_{-}^{*}\cE' \to \pi_{-}^{*}\cE'
\]
the object $C$ is such that $C|_{\cX} \in \langle \cA_{i'}(i'-l), \ldots, \cA_{i'}(i'-1) \rangle$.

We now claim that $\cA_{i}(i) \subseteq ^{\perp}\!\!\langle \cA_{i'}(-k\cL, i'-l+k), \ldots \cA_{i'}(-k\cL, i'-1+k)\rangle$ for all $k \in \ZZ$.
Suppose this is not the case, then one can find $\cE \in \cA_{i}(i)$ and $\cF \in \cA_{i'}(-k\cL,j)$ with $j \in [i'-l+k, i'-1+k]$ such that $\Hom(\cE,\cF) \not= 0$.
We must then have $j = i$, which gives $i' \ge i-k + 1 > i - k \ge i'-l \ge 0$, and hence $k \ge i-i' +  1 \ge 1$.
But then $\cA_{i'}(-k\cL) \subseteq \cA_{i-k}(-k\cL) \subseteq \cA_{i}^{\perp}$.

It follows that $\Hom(\cE, C|_{\cX}\otimes \cO(-k\cL,k)) = 0$ for all $k \in \ZZ$, and by Lemma \ref{thm:UsefulTrickVariant}, this gives $\Hom(\cE, C) = 0$.
Hence $\Hom(\pi_{-}^{*}\cE, i_{-}^{*}(i_{-})_{!}\pi_{-}^{*}\cE') = \Hom(\pi_{-}^{*}\cE, \pi_{-}^{*}\cE')$, and so since $\Hom(\cE, \cE'(i\cL,-i)) = 0$ for $i \ge 1$, then by Lemma \ref{thm:UsefulTrickVariant} it follows that $\Hom(\pi_{-}^{*}\cE, \pi_{-}^{*}\cE') = \Hom(\cE,\cE')$.
\end{proof}

For any category $\cS \in D^{b}(\cX)$, let $\cW_{0}(S) \subset D^{b}(\wt \cX)$ be the subcategory of those objects $\cE$ such that $i_{0}^{*}\cE \in \cS$.
We write $j_{+}$ for the inclusion $\wt \cX \setminus \cX \into \wt \cX$.
\begin{nlemma}
\label{thm:restrictionIdentification}
The functor $j_{+}^{*} \circ \pi_{-}^{*} \colon D^{b}(X)(0) \to D^{b}(\wt \cX \setminus \cX)$ is an equivalence. 
Let $F$ denote its inverse.
If $\cE \in \cW_{0}(\cS_{+})$, we have $F(j_{+}^{*}\cE) \in \langle \cA_{i}(i\cL) \rangle_{i \in I}$ for some set $I$ if and only if $\cE \in \cW_{0}(\langle \cA_{i}(i) \rangle_{i \in I})$.
\end{nlemma}
\begin{proof}
The functor is an equivalence because $\pi_{-} \colon \wt\cX \setminus \cX \to X$ is an isomorphism.
Using Lemma \ref{thm:SODOfTildeX}, if $\cE \in \cW_{0}(\cS_{+})$, then it may be written as an iterated extension of objects of the form $\pi_{-}^{*}\cE_{i}$ with $\cE_{i} \in \cA_{i}(i)$.
Since $\cL|_{\wt \cX \setminus \cX} \cong \cO(1)|_{\wt \cX \setminus \cX}$, we have $F(\pi_{-}^{*}\cA_{i}(i)) = \cA_{i}(i\cL)$.
Thus we get an expression of $F(j^{*}\cE)$ as an iterated extension of objects $F(\cE_{i}) \in \cA_{i}(i\cL)$.
Since the expression of both $\cE$ and $F\pi_{-}^{*}\cE$ as such iterated extensions are unique, and since $F(\cE_{i}) = 0 \Leftrightarrow \cE_{i} = 0$, the claim follows.
\end{proof}

Recall that $\cC_{L} = \cW(\cS_{0}) = \cW(\langle \cA_{0}, \ldots \cA_{l-1}((l-1)\cL))$.
\begin{nprop}
\label{thm:SODOfPlusWindow}
There is a semiorthogonal decomposition of $\cD_{L^{\perp}} \cong \cW(\cS_{+})$ into
\[
\cW(\cS_{+}) = \langle \cC_{L}, \Psi\cA_{l}(l), \ldots, \Psi\cA_{n}(n) \rangle.
\]
\end{nprop}
\begin{proof}
The restriction of the functor $\Psi$ to $D^{b}(X)(0)$ admits a right adjoint $(-)^{\CC^{*}}\circ(\pi_{-})_{*}(i_{-})^{*}$.
Hence the restriction of $\Psi$ to $\cA_{i}(i)$ admits a right adjoint for all $i$, and so the $\Psi\cA_{i}(i)$ are right admissible subcategories.

It's therefore enough to show that $\cC_{L} = \langle \Psi\cA_{l}(l), \ldots, \Psi\cA_{n}(n) \rangle^{\perp}$.
Let $\cE \in \cA_{i}(i)$ with $i \in [l,n]$ and $\cF \in \cW(\cS_{+})$.
By adjunction, $\Hom(\Psi\cE, \cF) = \Hom(\pi_{-}^{*}\cE, i_{-}^{*}\cF)$.
By Lemma \ref{thm:localCohomology} (applied to the case $L = 0$), we get $j_{+} \colon \wt \cX \setminus \cX \to \wt \cX$ and get 
\[
\Hom(\pi_{-}^{*}\cE, i_{-}^{*}\cF) = \Hom(j_{+}^{*}\pi_{-}^{*}\cE, j_{+}^{*}i_{-}^{*}\cF).
\]
Using the equivalence $F \colon D^{b}(\wt \cX \setminus \cX) \to D^{b}(X)(0)$ of Lemma \ref{thm:restrictionIdentification}, this Hom space equals $\Hom(\cE(i\cL,-i), Fj_{+}^{*}i_{-}^{*}\cF)$.
The claim that the Hom space vanishes for all $\cE \in \cA_{i}(i)$ with $i \in [l,n]$ is equivalent to the claim that 
\[
Fj_{+}^{*}i_{-}^{*}\cF \in \langle \cA_{l}(l\cL), \ldots, \cA_{n}(n\cL) \rangle^{\perp} = \langle \cA_{0}, \ldots, \cA_{l-1}((l-1)\cL)\rangle.
\]
By Lemma \ref{thm:restrictionIdentification}, this is equivalent to $\cF \in \cW(\cA_{0}, \ldots, \cA_{l-1}((l-1)\cL) = \cC_{L}$.
\end{proof}

\begin{nremark}
In \cite{kuznetsov_homological_2007}, the subcategories $\cA_{i}(i\cL) \subseteq D^{b}(X_{L^{\perp}})$ are defined as the images along the restriction functor $D^{b}(X) \to D^{b}(X_{L^{\perp}})$. 
One can check that when we have the Kn\"orrer equivalence $\cD_{L^{\perp}} \cong D^{b}(X_{L^{\perp}})$, the decomposition of $\cD_{L^{\perp}}$ that we find is the same as the one used in that paper.

The existence of the semi-orthogonal decomposition of $\cD_{L^{\perp}}$ could have been obtained more simply by working directly on $D^{b}(\cZ_{+},W)$ rather than in $\cW(\cS_{+})$; we've gone to some extra trouble in order to identify the piece $\cC_{L}$ with $\cW(\cS_{0})$.
\end{nremark}

\subsection{The equivalence $\cW(\cS_{-}) \to \cD_{L}^{\vee}$}
Let $\cC \subset D^{b}(X)$ be a full, triangulated subcategory, and define $D^{b}(\cZ_{-},W)_{\cC} \subseteq D^{b}(\cZ_{-},W)$ as the subcategory of those objects $\cE$ such that for all $[f] \in \PP L$, we have $\cE|_{[f] \times X} \in \cC$.
Note that this generalises our definition of $\cD_{L}^{\vee}$, because we have $D^{b}(\cZ_{-},W)_{\cA_{0}} = \cD_{L}^{\vee}$
\begin{nlemma}
\label{thm:fullyFaithfulnessZMinus}
If $\cC \subseteq D^{b}(X)$ is admissible, then the restriction functor $j_{-}^{*}$ gives an equivalence $\cW(\cC(0), \ldots, \cC(l-1)) \cong D^{b}(\cZ_{-},W)_{\cC}$.
\end{nlemma}
\begin{proof}
Using Prop.~\ref{thm:GITResultOurCase}, we immediately get fully faithfulness, and the fact that for $\cE \in D^{b}(\cZ_{-},W)$, there exists a unique $\cF \in \cW(D^{b}(X)(0), \ldots D^{b}(X)(l-1)\rangle$ such that $j_{+}^{*}\cF = \cE$.
It remains to see that $\cE \in D^{b}(\cZ_{-},W)_{\cC} \Leftrightarrow \cF \in \cW(\cC(0),\ldots,\cC(l-1))$.

So let $\cF \in \cW(D^{b}(X)(0), \ldots D^{b}(X)(l-1)\rangle$.
By Lemma \ref{thm:SODOfL}, we have $i_{-}^{*}\cF \in \langle \pi_{-}^{*}D^{b}(X)(l-1), \ldots, \pi_{-}^{*}D^{b}(X)(0) \rangle$.
The claim to show is that $i_{-}^{*}\cF \in \langle \pi_{-}^{*}\cC(l-1), \ldots, \pi_{-}^{*}\cC(0) \rangle$ if and only if $i_{-}^{*}\cF|_{[f] \times X} \in \cC$ for all $[f] \in \PP L$.

Since $\cC$ is admissible, we may refine the semi-orthogonal decomposition as follows:
\begin{align*}
& \langle \pi_{-}^{*}D^{b}(X)(l-1), \ldots, \pi_{-}^{*}D^{b}(X)(0) \rangle \\
=& \langle \pi_{-}^{*}\cC^{\perp}(l-1), \ldots, \pi_{-}^{*}\cC^{\perp}(0), \pi_{-}^{*}\cC(l-1), \ldots, \pi_{-}^{*}\cC(0) \rangle = \langle \cT_{\cC^{\perp}} , \cT_{\cC}\rangle,
\end{align*}
where $\cT_{\cC^{\perp}}$ (resp.~$\cT_{\cC}$) denotes the span of the first (resp.~last) $l$ pieces of the decomposition.

Let $(i_{-}^{*}\cF)_{\cT_{\cC}} \to i_{-}^{*}\cF \to (i_{-}^{*}\cF)_{\cT_{\cC^{\perp}}}$ be the corresponding decomposition of $i_{-}^{*}\cF$.
It's clear that for any $[f] \in \PP L$, we have $(i_{-}^{*}\cF)_{\cT_{\cC}}|_{[f] \times X} \in \cC$.
If $\cF \in \cW(\cS_{-})$, then $i_{-}^{*}\cF = (i_{-}^{*}\cF)_{\cT_{\cC}}$, hence $i_{-}^{*}\cF|_{[f] \times X} \in \cC$. 
If $\cF \not\in \cW(\cC)$, then $(i_{-}^{*}\cF)_{\cT_{\cC^{\perp}}} \not= 0$, hence there is some point $[f] \in \PP L$ such that $(i_{-}^{*}\cF)_{\cT_{\cC^{\perp}}}|_{[f] \times X} \not= 0$.
But since $(i_{-}^{*}\cF)_{\cT_{\cC^{\perp}}}|_{[f] \times X} \in \cC^{\perp}$, it follows that $(i_{-}^{*}\cF)|_{[f] \times X} \not\in \cC$. 
\end{proof}

\begin{ncor}
\label{thm:windowEquivalenceMinus}
The restriction functor $j_{-}^{*}$ gives an equivalence $\cW(\cS_{-}) \cong \cD_{L}^{\vee}$.
\end{ncor}

\subsection{The left Lefschetz decomposition of $\cD$}
\label{sec:dualLefschetz}
The only thing that remains in the proof of Thm.~\ref{thm:fundamentalTheoremOfHPD} is to describe the orthogonal complement to $\cC_{L}$ in $\cD_{L}^{\vee}$.
As a step towards this, we'll need to produce a new Lefschetz decomposition of $\cD$.
To be precise, the new decomposition is a left Lefschetz decomposition, as opposed to the standard right Lefschetz decomposition, which means that the bigger categories are left orthogonal rather than right orthogonal to the smaller categories.

Recall from \cite[Sec.~4]{kuznetsov_homological_2007} that there exists a semi-orthogonal decomposition $\cA_{0} = \langle \fa_{0}, \ldots, \fa_{n} \rangle$, defined by $\langle \fa_{i}, \ldots \fa_{n} \rangle = \cA_{i}$.
Kuznetsov also produces a second decomposition $\langle \alpha^{*}(\fa_{0}(\cL)), \ldots, \alpha^{*}(\fa_{n}(n\cL)) \rangle$, where $\alpha \colon \cA_{0} \to \cD$ is the inclusion functor and $\alpha^{*} \colon \cD \to \cA_{0}$ is the left adjoint \cite[Lem.~4.3]{kuznetsov_homological_2007}.
The proof given there uses a Serre functor, which we don't assume to exist on $\cD$, so let's give a modified, Serre-functor-less proof of the statement.
\begin{nlemma}
The functor $\alpha^{*} \colon \fa_{i}((i+1)\cL) \to \cA_{0}$ is fully faithful for all $i$, and this gives a semi-orthogonal decomposition
\[
\cA_{0} = \langle \alpha^{*}\fa_{0}(\cL), \ldots, \alpha^{*}\fa_{n}((n+1)\cL).
\] 
\end{nlemma}
\begin{proof}
Let $\cE_{i} \in \mathfrak{a}_{i}((i+1)\cL)$ and $\cE_{j} \in \mathfrak{a}_{j}((j+1)\cL)$, with $i \ge j$.
In the proof of \cite[Lem.~4.2]{kuznetsov_homological_2007}, it is shown that 
\begin{equation}
\label{eqn:factFromHPD}
\cE_{j} \in \langle \cA_{0}, \ldots, \cA_{j}(j\cL) \rangle \cap ^{\perp}\!\!\langle \cA_{1}(\cL), \ldots, \cA_{j}(j\cL) \rangle.
\end{equation}
Thus the cone $C$ over $\cE_{j} \to \alpha^{*}\cE_{j}$ is contained in $^{\perp}\cA_{0} \cap \langle \cA_{0}, \ldots, \cA_{j}(j\cL) \rangle = \langle \cA_{1}(\cL), \ldots, \cA_{j}(j\cL) \rangle$.
By \eqref{eqn:factFromHPD} with $i$ instead of $j$ we have $\cE_{i} \in ^{\perp}\!\!\langle \cA_{1}(\cL), \ldots, \cA_{i}(i\cL) \rangle$, and so $\Hom(\cE_{i},C) = 0$.
Hence $\Hom(\alpha^{*}\cE_{i}, \alpha^{*}\cE_{j}) = \Hom(\cE_{i}, \alpha^{*}\cE_{j}) = \Hom(\cE_{i}, \cE_{j})$.

Taking $i = j$, this shows that $\alpha^{*} \colon \fa_{i}((i+1)\cL) \to \cA_{0}$ is fully faithful, while taking $i > j$ shows the semi-orthogonality of the $\alpha^{*}(\fa_{i}((i + 1)\cL))$.
The fact that these categories generate $\cA_{0}$ is shown as in \cite[Lem.~4.3]{kuznetsov_homological_2007}.
\end{proof}

\begin{nlemma}\cite[Lem.~4.5]{kuznetsov_homological_2007}
\label{thm:ANoughtsGenerate}
We have $\langle \cA_{0}, \ldots, \cA_{i}(i\cL) \rangle = \langle \cA_{0}, \ldots, \cA_{0}(i\cL) \rangle$.
\end{nlemma}

Define the category 
\[
\cA_{i}^{\dagger} = \cA_{0} \cap \langle \cA_{0}(\cL), \ldots, \cA_{0}((i-1)\cL) \rangle^{\perp} = \langle \alpha^{*}(\fa_{0}(\cL)), \ldots, \alpha^{*}(\fa_{i-1}((i-1)\cL)) \rangle^{\perp}.
\]
\begin{nlemma}
There is a semiorthogonal decomposition $\cD = \langle \cA_{n}^{\dagger}(-n\cL), \ldots, \cA_{0}^{\dagger} \rangle$.

For any $i \le n$, we have $\langle \cA_{i}^{\dagger}(-i\cL), \ldots, \cA_{0}^{\dagger} \rangle = \langle \cA_{0}(-i\cL), \ldots, \cA_{0} \rangle$.
\end{nlemma}
\begin{proof}
Semi-orthogonality is straightforward, so we only prove generation.

We prove by induction that $\langle \cA_{i}^{\dagger}(-i\cL), \ldots, \cA_{0}^{\dagger} \rangle = \langle \cA_{0}(-i\cL), \ldots, \cA_{0} \rangle$.
Let $\cE \in \alpha^{*}(\fa_{j-1}(j\cL))$, and let $\cE' \in \fa_{j-1}(j\cL)$ be such that $\alpha^{*}\cE' = \cE$.
In the exact triangle $C \to \cE' \to \cE$, we have 
\begin{align*}
C \in ^{\perp}\!\!\cA_{0} \cap \langle \cA_{0}, \cA_{0}(j\cL) \rangle &\subseteq ^{\perp}\!\!\cA_{0} \cap \langle \cA_{0}, \cA_{1}(\cL), \ldots, \cA_{j}(j\cL) \rangle \\
& = \langle \cA_{1}(\cL), \ldots, \cA_{j}(j\cL) \rangle,
\end{align*}
using Lemma \ref{thm:ANoughtsGenerate}. 
Hence $\cE \in \langle \cA_{0}(\cL), \ldots, \cA_{0}(j\cL) \rangle$.

It follows that for $j \le i$, we have $\alpha^{*}(\fa_{j-1}(j\cL))(-i\cL) \subseteq \langle \cA_{0}(-(i-1)\cL), \ldots, \cA_{0}\rangle$.
We then find 
\begin{align*}
\langle \cA_{i}^{\dagger}(-i\cL), \ldots, \cA_{0}^{\dagger} \rangle &= \langle \cA_{i}^{\dagger}(-i\cL), \cA_{0}(-(i-1)\cL), \ldots, \cA_{0} \rangle \\
&= \langle \cA_{0}(-i\cL), \cA_{0}(-(i-1)\cL), \ldots, \cA_{0} \rangle,
\end{align*}
using induction and the fact that $\cA_{0} = \langle \cA_{i}^{\dagger}, \alpha^{*}(\fa_{0}(\cL)), \ldots, \alpha^{*}(\fa_{i-1}(i\cL))\rangle$.
\end{proof}

\subsection{The full decomposition of $\cD^{\vee}_{L}$}
The only thing that remains to prove in Thm.~\ref{thm:fundamentalTheoremOfHPD} is the existence of the semi-orthogonal pieces $\cB_{i}(i)$ in the decomposition of $\cD_{L}^{\vee}$.
In this section, we'll assume that $l = \dim L \ge 2$, since otherwise $\cD_{L}^{\vee} = \cC_{L}$.

We let $\cA_{i}^{*} = \cA_{0} \cap \cA_{i}^{\perp}$.

\begin{nlemma}
\label{thm:LimitingTheWeightsOfAStar}
For $i = 1, \ldots, n$, we have 
\[
\cA_{i}^{*}(\cL) \subseteq  \langle \cA_{i-1}^{\dagger}(-(i-1)\cL), \ldots, \cA_{0}^{\dagger} \rangle.
\]
\end{nlemma}
\begin{proof}
We have $\cA_{i}^{*} \subseteq \cA_{0} \subseteq \langle \cA_{i+1}(\cL), \ldots, \cA_{n}((n-i)\cL)\rangle^{\perp}$, and also $\cA_{i}^{*} \subseteq \cA_{i}^{\perp}$.
Hence
\begin{align*}
\cA_{i}^{*}(\cL) \subseteq \langle \cA_{i}(\cL), \cA_{i+1}(2\cL), \ldots, \cA_{n}((n-i+1)\cL)\rangle^{\perp} &= \langle \cA_{0}(-(i-1)\cL), \ldots, \cA_{i-1} \rangle \\
&= \langle \cA_{0}(-(i-1)\cL), \ldots, \cA_{0} \rangle \\
&= \langle \cA_{i-1}^{\dagger}(-(i-1)\cL), \ldots, \cA_{0}^{\dagger} \rangle.
\end{align*}
\end{proof}

Let $\Omega = (i_{+})_{*}\pi_{+}^{*} \colon D^{b}(\cX) \to D^{b}(\cZ,W)$, and let $\Omega_{-} = j_{-}^{*}(i_{+})_{*}\pi_{+}^{*} \colon D^{b}(\cX) \to D^{b}(\cZ_{-},W)$.
The functor $\Omega$ restricted to $D^{b}(X)(0)$ admits a left adjoint given by
\begin{align*}
(-)_{X}^{\vee} \circ (-)^{\CC^{*}} \circ (\pi_{+})_{*} \circ (-)_{\cY_{+}}^{\vee} \colon D^{b}(\cY_{+}) \to D^{b}(X).
\end{align*}
The next three lemmas are proved by the same arguments as Lemmas \ref{thm:InclusionIntoZMinusFullyFaithful}, \ref{thm:SODOfZ} and \ref{thm:projectionOfWeights}.
The left adjoint to $\Omega|_{D^{b}(X)(0)}$ above is used in order to get left admissibility of $\Omega(\cA^{\dagger}_{i}(-i\cL,j))$ in Lemma \ref{thm:SODOfZDual}.

\begin{nlemma}
\label{thm:OmegaFullyFaithful}
The functor $\Omega \colon D^{b}(X)(k) \to D^{b}(\cZ,W)$ is fully faithful for any $k \in \ZZ$.  
\end{nlemma}

\begin{nlemma}
\label{thm:SODOfZDual}
There exists a semi-orthogonal decomposition
\[
\cW(\langle \cD(i) \rangle_{i \in \ZZ}) \cong \langle \langle \Omega(\cA^{\dagger}_{i}(-i\cL,j)) \rangle_{i,j}, \cW(\langle\cA_{0}(i)\rangle_{i \in \ZZ})  \rangle,
\]
where $i = 1, \ldots, n$ and $j \in \ZZ$.
The semi-orthogonality relations are such that $\Omega(\cA^{\dagger}_{i}(-i\cL,j)) \subseteq \Omega(\cA^{\dagger}_{i'}(-i'\cL,j'))^{\perp}$ if
\begin{itemize}
\item $i > i'$, or if
\item $i = i'$ and $j > j'$.
\end{itemize}
\end{nlemma}
Recall that $\iota \colon \cW(\langle \cA_{0}(i)\rangle_{i\in\ZZ}) \to D^{b}(\cZ,W)$ is the inclusion functor.
\begin{nlemma}
\label{thm:projectionOfWeightsZDual}
Let $\cE \in \cW(\langle \cA_{i}^{\dagger}(-i\cL,j) \rangle_{(i,j) \in I}$, where $I \subseteq [0,n] \times \ZZ$.
Then $\iota^{!}\cE \in \cW(\langle \cA_{i}(j-i)\rangle_{(i,j) \in I})$.

Let $I' = I \cap ([1,n] \times \ZZ)$, and let $K = \bigcup_{(i,j) \in I'} [j-i, j] \cap \ZZ$.
The cone over $\iota^{!}\cE \to \cE$ lies in $\cW(\cD(k))_{k \in K}$.
\end{nlemma}

\begin{nlemma}
\label{thm:OmegaMinusFullyFaithful}
The functor $\Omega_{-} \colon D^{b}(X)(k) \to D^{b}(\cZ_{-},W)$ is fully faithful for any $k \in \ZZ$.
\end{nlemma}
\begin{proof}
We may assume $k = 0$, so let $\cE, \cF \in D^{b}(X)(0)$.
Then it's easy to see that $\Omega\cE, \Omega\cF \in \cW(D^{b}(X)(-1), D^{b}(X))$, and from this it follows that $\Hom(i_{0}^{*}\Omega \cE, i_{0}^{*}\Omega \cF(i)) = 0$ for $i \le -2$, and thus for $i \le -l$.
Hence by Lemmas \ref{thm:bestLocalCohomologyStatement} and \ref{thm:OmegaFullyFaithful}, we have $\Hom(\Omega \cE, \Omega \cF) = \Hom(j_{-}^{*}\Omega \cE, j_{-}^{*}\Omega \cF) =  \Hom(\Omega_{-}\cE, \Omega_{-}\cF)$.
\end{proof}

\begin{nlemma}
\label{thm:SODOfZMinus}
There exists a semi-orthogonal decomposition
\[
D^{b}(\cZ_{-},W)_{\cD} \cong \langle \langle \Omega_{-}(\cA^{\dagger}_{i}(-i\cL,j)) \rangle_{i,j}, \cD^{\vee}_{L}  \rangle,
\]
where $i = 1, \ldots, n$ and $j = 0, \ldots, l-1$.
The semi-orthogonality relations are such that $\Omega_{-}(\cA^{\dagger}_{i}(-i\cL,j)) \subseteq \Omega_{-}(\cA^{\dagger}_{i'}(-i'\cL,j'))^{\perp}$ if
\begin{itemize}
\item $i > i'$, or if
\item $i = i'$ and $j > j'$.
\end{itemize}
\end{nlemma}
\begin{proof}
Possibly twisting the window subcategory $\cW(\cS_{-})$ by $\cO(-1)$, by Prop.~\ref{thm:OurSpecialGITStatement} we find that $j_{+}^{*}$ admits a partial inverse $(j_{+}^{*})^{-1} \colon D^{b}(\cZ,W) \to \cW(\cS_{-})$, moreover this can be chosen such that $(j_{-}^{*})^{-1}\Omega_{-} = \Omega \colon \cA_{i}^{\dagger}(-i\cL,j) \to D^{b}(\cZ,W)$.
This implies that $\Omega_{-}$ admits a left adjoint when restricted to $\cA^{\dagger}_{i}(-i\cL,j)$, and so the $\Omega_{-}(\cA^{\dagger}_{i}(-i\cL,j))$ are left admissible.

Let $\cE \in \cA^{\dagger}_{i}(-i\cL,j)$, $\cF \in \cA^{\dagger}_{i'}(-i'\cL,j')$, and assume that either of the two conditions listed hold.
By Lemma \ref{thm:SODOfZDual}, we have $\Hom(\Omega\cF, \Omega\cE)=0$.

We have 
\[
\Omega\cE \in \cW(\cA^{\dagger}_{i}(-i\cL,j), \cA^{\dagger}_{i}(-(i-1)\cL,j-1))
\]
and 
\[
\Omega\cF \in \cW(\cA^{\dagger}_{i'}(-i'\cL,j'), \cA^{\dagger}_{i'}(-(i'-1)\cL,j'-1)).
\]
It's easy to check that $\langle \cA^{\dagger}_{i}(-i\cL,j), \cA^{\dagger}_{i}(-(i-1)\cL,j-1) \rangle \subseteq \langle \cA^{\dagger}_{i'}(-i'\cL,j' + k), \cA^{\dagger}_{i'}(-(i'-1)\cL,j'-1 +k) \rangle^{\perp}$ for $k \le -l$, and so $\Hom(i_{0}^{*}\Omega \cF, i_{0}^{*}\Omega \cE(k)) = 0$ for $k \le -l$.
Hence by Lemma \ref{thm:bestLocalCohomologyStatement}, $\Hom(\Omega_{-}\cF,\Omega_{-}\cE) = \Hom(\Omega\cF,\Omega\cE) = 0$.

Now let $\cF \in D^{b}(\cZ_{-},W)_{\cD}$, and let $\cF' \in \cW(\cD(0), \ldots, \cD(l))$ be such that $\cF = j_{-}^{*}\cF'$.
If $\cE \in \cA^{\dagger}_{i}(-i\cL,j)$, then using Lemma \ref{thm:bestLocalCohomologyStatement}, we find $\Hom(\cF, \Omega_{-}\cE) = \Hom(\cF', \Omega\cE)$.
Hence if this Hom space vanishes for all such $\cE \in \cA_{i}^{\dagger}(-i\cL,j)$, Lemma \ref{thm:SODOfZDual} shows that $\cF' \in \cW(\cA_{0}(0),\ldots, \cA_{0}(l-1))$, and so $\cF \in \cD_{L}^{\vee}$.
\end{proof}

Let $\iota_{-} \colon \cD_{L}^{\vee} \to D^{b}(\cZ_{-},W)_{\cD}$ denote the inclusion.
\begin{nlemma}
\label{thm:RightAdjointsCommute}
If $\cE \in \cA_{i}^{*}(i)$, then $\iota_{-}^{!}\Omega_{-}\cE = j_{-}^{*}\iota^{!}\Omega\cE$.
\end{nlemma}
\begin{proof}
By Lemmas \ref{thm:LimitingTheWeightsOfAStar} and \ref{thm:projectionOfWeightsZDual}, we have $\iota^{!}\Omega \cE \in \cW(\cA_{0}, \ldots, \cA_{0}(i))$.
Let now $\cF \in \cD_{L}^{\vee}$, and let $\cF' \in \cW(\cA_{0}, \ldots ,\cA_{0}(l-1))$ be such that $j^{*}_{-}\cF' = \cF$.
Hence by Lemma \ref{thm:bestLocalCohomologyStatement}, we have 
\begin{align*}
\Hom(\cF, j_{-}^{*}\iota^{!}\Omega\cE) &= \Hom(j_{-}^{*}\cF', j_{-}^{*}\iota^{!}\Omega\cE) = \Hom(\cF', \iota^{!}\Omega\cE) = \Hom(\cF', \Omega\cE) \\
&= \Hom(\cF, \Omega_{-}\cE) = \Hom(\cF,\iota_{-}^{!}\Omega_{-}\cE),
\end{align*}
and so $\iota_{-}^{!}\Omega_{-}\cE = j_{-}^{*}\iota^{!}\Omega\cE$.
\end{proof}

\begin{nlemma}
\label{thm:projectionOfWeightsZMinus}
\label{thm:ImageOfBProjectionDual}
If $\cE \in \cA_{i}^{*}(i)$, then $\iota_{-}^{!}\Omega_{-}\cE \in j_{-}^{*}(\cW(\cA_{0}, \ldots, \cA_{0}(i)))$, and the cone over $\Omega_{-}\cE \to \iota_{-}^{!}\Omega_{-}\cE$ lies in $j_{-}^{*}\cW(\cD, \ldots, \cD(i-1))$.
\end{nlemma}
\begin{proof}
Combine Lemmas \ref{thm:LimitingTheWeightsOfAStar}, \ref{thm:projectionOfWeightsZDual} and \ref{thm:RightAdjointsCommute}.
\end{proof}

\begin{nlemma}
\label{thm:FullyFaithfulDualCollection}
The functor $\iota_{-}^{!}\Omega_{-} \colon \cA^{*}_{i}(i) \to D^{b}(\cZ_{-},W)$ is fully faithful.
If $i \le i'$, then $\iota_{-}^{!}\Omega_{-}(\cA^{*}_{i}(i)) \subseteq ^{\perp}\!\! \iota^{!}\Phi(\cA^{*}_{i'}(i'))$.
\end{nlemma}
\begin{proof}
Let $\cE, \cE' \in \cA_{i}^{*}(i)$.
By adjunction, $\Hom(\iota_{-}^{!}\Omega_{-}\cE, \iota_{-}^{!}\Omega_{-}\cE') = \Hom(\iota_{-}^{!}\Omega_{-}\cE, \Omega_{-}\cE')$.
By Lemma \ref{thm:ImageOfBProjectionDual}, the cone $C$ over $\iota_{-}^{!}\Omega_{-}\cE \to \Omega_{-}\cE$ lies in $j_{-}^{*}(\cW(\cD(0), \ldots, \cD(i-1)))$.
Hence $\Hom(C, \Omega\cE') = \Hom(i_{+}^{*}C, \pi_{+}^{*}\cE') = 0$, using Lemmas \ref{thm:HelpfulTrick} and \ref{thm:bestLocalCohomologyStatement}.
It follows that $\Hom(\iota_{-}^{!}\Omega_{-}\cE, \Omega_{-}\cE') = \Hom(\Omega_{-}\cE, \Omega_{-}\cE') = \Hom(\cE,\cE')$, by Lemma \ref{thm:OmegaMinusFullyFaithful}.

A similar and easier argument shows that if $\cE \in \cA^{*}_{i}(i)$ and $\cE' \in \cA^{*}_{i'}(i')$ with $i < i'$, then $\Hom(\iota_{-}^{!}\Omega_{-}\cE, \iota_{-}^{*}\Omega_{-}\cE') = 0$.
\end{proof}

We define the subcategory $\cB_{i} \subset \cD^{\vee}_{L}$ as the image of $\iota_{-}^{!}\Omega_{-}\cA^{*}_{i}$.
Note that the functors $\Omega_{-}$ and $\iota_{-}^{!}$ both commute with $- \otimes \cO(1)$, so that $\cB_{i}(i) = \iota_{-}^{!}\Omega_{-}\cA^{*}_{i}(i)$.

\begin{nprop}
There exists a semiorthogonal decomposition
\[
\cD^{\vee}_{L} = \langle \cB_{l-1}(l-1), \ldots, \cB_{1}(1), j_{-}^{*}\cC_{L} \rangle.
\]
\end{nprop}
\begin{proof}
By Lemma \ref{thm:FullyFaithfulDualCollection}, the categories $\cB_{i}(i)$ are semi-orthogonal to each other.
Since both $\Omega_{-}$ and $\iota_{-}^{!}$ admit left adjoints, the $\cB_{i}$ are left admissible, and it suffices to show that $j_{-}^{*}\cC_{L} = ^{\perp}\!\!\langle \cB_{l}(l), \ldots, \cB_{l}(1) \rangle$.

Let now $\cE \in \cD_{L}^{\vee}$ and $\cF \in \cA_{i}^{*}(i)$ with $i \ge 1$.
Then $\Hom(\cE, \iota_{-}^{!}\Omega_{-}\cF) = \Hom(\cE, \Omega_{-}\cF)$.
Let $\cE' \in \cW(\cS_{-})$ be such that $j_{-}^{*}\cE' = \cE$.
We then have 
\[
\Hom(\cE, \Omega_{-}\cF) = \Hom(j_{-}^{*}\cE', j_{-}^{*}\Omega\cF) = \Hom(\cE', \Omega\cF) = \Hom(i_{+}^{*}\cE', \pi_{+}^{*}\cF),
\]
using Lemma \ref{thm:bestLocalCohomologyStatement}.

Now if $\cE' \in \cC_{L}$, then $\Hom(i_{0}^{*}\cE', \cF(k)) = 0$ for all $k \ge 0$, and since this holds for any $\cF \in \cA_{i}^{*}(i)$, it follows that $\cE \in ^{\perp}\!\!\langle \cB_{l-1}(l-1), \ldots, \cB_{1}(1) \rangle$.

If $\cE' \not\in \cC_{L}$, let $j$ be maximal such that $i_{0}^{*}\cE'|_{\cD(j)} \not\in \cA_{j}(j)$, and let $\cF$ be the projection of $i_{0}^{*}\cE'|_{\cD(j)}$ to $\cA_{j}^{*}(j)$.
There is then by construction a non-trivial map in $\Hom(i_{0}^{*}\cE', \cF)$, and by the maximality of $j$ we find $\Hom(i_{0}^{*}\cE', \cF(k)) = 0$ for $k > 0$.
Hence $\Hom(\cE, \iota_{-}^{!}\Omega_{-}\cF) \not= 0$, and so $\cE \not\in ^{\perp}\!\!\langle \cB_{l}(l), \ldots, \cB_{1}(1) \rangle$.
\end{proof}

\section{Different notions of base change}
\label{sec:baseChange}
In this section, we assume that $X$ and $Y$ are smooth, projective varieties and $\cL$ is globally generated.
The result of the section is essentially that two natural notions of base change of a category agree, and it surely holds with weaker assumptions than these.

Let $\cH \subset X \times \PP V^{\vee}$ be the incidence variety of pairs $(x,H)$ with $x \in H$.
Using Prop.~\ref{thm:KnorrerPeriodicity}, we find that $D^{b}(\cH) \cong D^{b}(\cZ_{-},W)$, where $\cZ_{-} = [\wt X \times (V^{\vee}\setminus 0)/\CC^{*}]$ and $W$ is as before.
We may therefore think of $\cD^{\vee}$ as a subcategory of $D^{b}(\cH)$.
More generally, we can describe $\cD^{\vee}_{L}$ as a subcategory of $D^{b}(\cH_{L})$, where $\cH_{L} = \PP L \times_{\PP V^{\vee}} \cH$.

For each hyperplane $H \in \PP V^{\vee}$, let $i_{H} \colon X \cap H \into (X \cap H) \times [H] \subset \cH$ and $j_{H} \colon X \cap H \into X$ be the natural inclusions.
Then $\cD^{\vee}$ may be described as the subcategory of those $\cE$ such that for any point $H \in \PP V^{\vee}$, we have $(j_{H})_{*}i_{H}^{*}\cE \in \cA_{0}$.
Similarly $\cD^{\vee}_{L} \subseteq D^{b}(\cH_{L})$ is the subcategory of objects such that $(j_{H})_{*}i_{H}^{*}\cE \in \cA_{0}$ for all $H \in \PP L$.

We say that a variety $Y$ over $\PP V^{\vee}$ represents $\cD^{\vee}$, if there is a kernel object $\cK \in D^{b}(Y \times_{\PP V^{\vee}} \cH)$ inducing an equivalence $q_{*}(\cK \otimes p^{*}(-)) \colon D^{b}(Y) \to \cD^{\vee}$, where $p$ and $q$ are the projections from $Y \times_{\PP V^{\vee}} \cH$ to $Y$ and $\cH$, respectively.
\begin{nprop}
If $Y$ represents $\cD^{\vee}$ and $L \subseteq \PP V^{\vee}$ is such that $Y_{L}$ has the expected dimension, then $D^{b}(Y_{L}) \cong \cD^{\vee}_{L}$.
\end{nprop}
\begin{proof}
Kuznetsov defines in \cite{kuznetsov_base_2011} a general notion of base change of subcategories, and using \cite[Cor.~5.7]{kuznetsov_base_2011}, it's easy to see that $\cD^{\vee}_{L}$ is the base change of $\cD^{\vee}$ along $\PP L \to \PP V$.
In the terminology of \cite{kuznetsov_hyperplane_2006}, the condition that $Y_{L}$ has the expected dimension ensures that $\PP L \to \PP V^{\vee}$ is faithful for $Y \to \PP V^{\vee}$, by \cite[Cor.~2.27]{kuznetsov_hyperplane_2006}.
The claim that $D^{b}(Y_{L}) \cong \cD_{L}^{\vee}$ is then a consequence of \cite[Thm.~6.4]{kuznetsov_base_2011}.
\end{proof}

\bibliographystyle{alphaabbrv}

\bibliography{bibliography}

\newcommand{\etalchar}[1]{$^{#1}$}
\begin{thebibliography}{BDF{\etalchar{+}}13}

\bibitem[ADS15]{addington_pfaffian-grassmannian_2015}
N.~Addington, W.~Donovan, and E.~Segal.
\newblock The {Pfaffian}--{Grassmannian} equivalence revisited.
\newblock {\em Algebr. Geom.}, 2(3):332--364, 2015.

\bibitem[AG15]{arinkin_singular_2015}
D.~Arinkin and D.~Gaitsgory.
\newblock Singular support of coherent sheaves and the geometric {Langlands}
  conjecture.
\newblock {\em Sel. Math. New Ser.}, 21(1):1--199, 2015.

\bibitem[BDF{\etalchar{+}}13]{ballard_homological_2013}
M.~Ballard, D.~Deliu, D.~Favero, M.~U. Isik, and L.~Katzarkov.
\newblock Homological {Projective} {Duality} via {Variation} of {Geometric}
  {Invariant} {Theory} {Quotients}.
\newblock {\em arXiv:1306.3957 [math]}, 2013.

\bibitem[BFK12]{ballard_variation_2012}
M.~Ballard, D.~Favero, and L.~Katzarkov.
\newblock Variation of geometric invariant theory quotients and derived
  categories.
\newblock {\em arXiv:1203.6643 [math]}, 2012.

\bibitem[BVdB03]{bondal_generators_2003}
A.~Bondal and M.~Van~den Bergh.
\newblock Generators and representability of functors in commutative and
  noncommutative geometry.
\newblock {\em Mosc. Math. J.}, 3(1):1--36, 258, 2003.

\bibitem[CT15]{carocci_homological_2015}
F.~Carocci and Z.~Turcinovic.
\newblock Homological projective duality for linear systems with base locus.
\newblock {\em arXiv:1511.09398 [math]}, 2015.

\bibitem[EP15]{efimov_coherent_2015}
A.~I. Efimov and L.~Positselski.
\newblock Coherent analogues of matrix factorizations and relative singularity
  categories.
\newblock {\em Algebra Number Theory}, 9(5):1159--1292, 2015.

\bibitem[Hir16]{hirano_derived_2016}
Y.~Hirano.
\newblock Derived {Kn\"orrer} periodicity and {Orlov}'s theorem for gauged
  {Landau}-{Ginzburg} models.
\newblock {\em arXiv:1602.04769 [math]}, 2016.

\bibitem[HL15]{halpern-leistner_derived_2015}
D.~Halpern-Leistner.
\newblock The derived category of a {GIT} quotient.
\newblock {\em J. Amer. Math. Soc.}, 28(3):871--912, 2015.

\bibitem[Isi13]{isik_equivalence_2013}
M.~U. Isik.
\newblock Equivalence of the {Derived} {Category} of a {Variety} with a
  {Singularity} {Category}.
\newblock {\em Int Math Res Notices}, 2013(12):2787--2808, 2013.

\bibitem[Kuz06]{kuznetsov_hyperplane_2006}
A.~G. Kuznetsov.
\newblock Hyperplane sections and derived categories.
\newblock {\em Izv. Ross. Akad. Nauk Ser. Mat.}, 70(3):23--128, 2006.

\bibitem[Kuz07]{kuznetsov_homological_2007}
A.~Kuznetsov.
\newblock Homological projective duality.
\newblock {\em Publications mathématiques de l'IHÉS}, 105(1):157--220, 2007.

\bibitem[Kuz11]{kuznetsov_base_2011}
A.~Kuznetsov.
\newblock Base change for semiorthogonal decompositions.
\newblock {\em Compos. Math.}, 147(3):852--876, 2011.

\bibitem[Kuz14]{kuznetsov_semiorthogonal_2014}
A.~Kuznetsov.
\newblock Semiorthogonal decompositions in algebraic geometry.
\newblock {\em arXiv:1404.3143 [math]}, 2014.

\bibitem[Orl12]{orlov_matrix_2012}
D.~Orlov.
\newblock Matrix factorizations for nonaffine {LG}-models.
\newblock {\em Math. Ann.}, 353(1):95--108, 2012.

\bibitem[Pos11]{positselski_two_2011}
L.~Positselski.
\newblock Two kinds of derived categories, {Koszul} duality, and
  comodule-contramodule correspondence.
\newblock {\em Mem. Amer. Math. Soc.}, 212(996):vi+133, 2011.

\bibitem[Ren15]{rennemo_homological_2015}
J.~V. Rennemo.
\newblock The homological projective dual of $\mathrm{Sym}^2 \mathbb{P}({V})$.
\newblock {\em arXiv:1509.04107 [math]}, 2015.

\bibitem[RS16]{rennemo_hori-mological_2016}
J.~V. Rennemo and E.~Segal.
\newblock Hori-mological projective duality.
\newblock {\em arXiv:1609.04045 {[hep-th]}}, 2016.

\bibitem[Seg11]{segal_equivalence_2011}
E.~Segal.
\newblock Equivalence between {GIT} quotients of {Landau}-{Ginzburg}
  {B}-models.
\newblock {\em Comm. Math. Phys.}, 304(2):411--432, 2011.

\bibitem[Shi12]{shipman_geometric_2012}
I.~Shipman.
\newblock A geometric approach to {Orlov}'s theorem.
\newblock {\em Compos. Math.}, 148(5):1365--1389, 2012.

\bibitem[Tel00]{teleman_quantization_2000}
C.~Teleman.
\newblock The quantization conjecture revisited.
\newblock {\em Ann. of Math.}, 152(1):1--43, 2000.

\bibitem[Tho15]{thomas_notes_2015}
R.~P. Thomas.
\newblock Notes on {HPD}.
\newblock {\em arXiv:1512.08985 [math]}, 2015.

\end{thebibliography}

\end{document}